\documentclass{elsart}
\usepackage{amsmath,amsfonts,amssymb}

\def\R{\mathbb R}
\def\C{\mathbb C}
\def\N{\mathbb N}
\def\Z{\mathbb Z}
\def\Sc{\mathcal S}
\def\H{\mathcal H}
\def\S{\mathbb S}
\def\d{\mathrm d}
\def\D{\mathrm D}

\def\br#1{{\left[{#1}\right]}}
\newcommand\Rn{{{\mathbb R}^n}}
\newcommand\pa{\partial}
\def\abs#1{{\left|{#1}\right|}}
\def\p#1{{\left({#1}\right)}}
\def\Fcal{\mathcal{F}}
\newcommand\Snm{{{\mathbb S}^{n-1}}}
\def\va{\varphi}
\def\i{\mathrm i}  
\let\Im\relax

\DeclareMathOperator{\Im}{Im}
\DeclareMathOperator{\diag}{diag}
\DeclareMathOperator{\supp}{supp}

\DeclareMathOperator{\spec}{spec}

\begin{document}

\begin{frontmatter}
\title{Dispersive estimates 
for hyperbolic systems with time-dependent coefficients}
\author[IC]{Michael Ruzhansky\thanksref{th:grant1}}
\author[ST]{Jens Wirth\thanksref{th:grant2}}
\thanks[th:grant1]{The first
 author was supported in part by the EPSRC
 grants EP/G007233/1 and EP/E062873/1.}
 \thanks[th:grant2]{The second author
 was supported by the EPSRC grant EP/E062873/1.}
\address[IC]{Department of Mathematics,\\
Imperial College London,\\
180 Queen's Gate, London SW7 2AZ, UK\\
m.ruzhansky@imperial.ac.uk}
\address[ST]{Institut f\"ur Analysis, Dynamik und Modellierung,\\Fachbereich Mathematik, Universit\"at Stuttgart,\\
Pfaffenwaldring 57, 70569 Stuttgart, Germany\\
jens.wirth@iadm.uni-stuttgart.de}

\begin{abstract}
This paper is devoted to the study of time-dependent 
hyperbolic systems and
the derivation of dispersive estimates for their solutions. 
It is based on a diagonalisation of the full symbol within 
adapted symbol classes in order to extract the essential 
information about representations of solutions. This is 
combined with a multi-dimensional van der Corput lemma 
to derive dispersive estimates.
\end{abstract}
\begin{keyword}
hyperbolic systems, 
diagonalisation, dispersive estimates, oscillating 
time-dependent coefficients, van der Corput lemma
\MSC 35L05\sep  35L15
\end{keyword}
\end{frontmatter}

\section{Introduction}

Our aim is to consider Cauchy problems for hyperbolic systems
\begin{equation}\label{eq:CP}
\D_t U =A(t,\D) U,\qquad  U(0,\cdot)=U_0,
\end{equation}
where $A(t,\D)$ denotes a smoothly time-dependent matrix Fourier multiplier with first order
symbol 
\begin{equation}
   A(t,\xi)\in C^\infty(\R_+\times\R^n_\xi,\C^{m\times m}) \cap L^\infty(\R_+,S^1(\R^n)) 
\end{equation}
subject to certain (natural) assumptions which are  described later on in detail. As usual we denote $\D_t=-\i\partial_t$
and we assume that eigenvalues of $A$ are real modulo `lower order terms'.

The decay of solutions to the 
scalar higher order equations with constant coefficients
of dissipative type was considered in \cite{RS05}.
Consequently, a comprehensive analysis of the dispersive
and Strichartz estimates of scalar hyperbolic equations
of higher orders with constant coefficients,
of general form, as well as for
hyperbolic systems with constant coefficients, 
with applications to nonlinear
equations was carried out in \cite{RS10}. In analogy to
this, the results of this paper have the corresponding
(by now standard) consequences for Strichartz estimates, 
and for the global
in time well-posedness of nonlinear equations, which we
omit here, and refer to \cite{RS10} for the corresponding
constructions.

At the same time, equations of the second order with
time dependent coefficients have been intensively
studied in the literature, see e.g.
\cite{RS05a}, 
\cite{Reissig2010}, \cite{RY00},
\cite{Wir06}, \cite{Wir10a}.
In \cite{MR10}, the asymptotic integration method was
developed for scalar hyperbolic equations with time
dependent homogeneous
symbols. The emphasis there was placed on minimising the
time-differentiability assumptions on the symbol, required
by applications to the Kirchhoff equations
(\cite{MR09}). Indeed, the
representation of solutions and dispersive estimates  
have been derived in \cite{MR10} for equations of any order,
by making assumptions on one 
time-derivative of the coefficients only, see
Remark \ref{REM:MR}. However, the asymptotic integration
produces an additional loss of regularity for high 
frequencies due to the fact that the amplitudes in the
representation of solutions are symbols of type $(0,0)$.
This is avoided in this paper since we construct the
amplitudes as symbols of type $(1,0)$.
In \cite{dAbbico:2009}
a special case of the results of this paper concerning differential hyperbolic 
systems with $\gamma=2$ (i.e., with non-vanishing Gaussian curvature of the Fresnel surfaces) 
has been treated.

Our approach is based on diagonalising the (full) symbol 
of the operator, modulo controllable terms, in order to gain a 
representation of solutions in terms of Fourier integrals. The basic
idea of this diagonalisation scheme comes from the treatment of
degenerate hyperbolic problems suggested in thee book of
Yagdjian \cite{Yagdjian}.
This representation is given in Theorem \ref{thm:repSol}.
Consequently, it is used to derive the dispersive estimates
for solutions to the Cauchy problem
\eqref{eq:CP} in Theorem \ref{THM:main} and
Theorem \ref{THM:main-nonconvex}, which are the main
results of this paper. For this purpose, we derive more
general estimates for oscillatory integrals in
Propositions \ref{PROP:oscillatory-convex}
and \ref{PROP:oscillatory-nonconvex}. 
The time decay rate of solutions to \eqref{eq:CP} is related to
the geometric properties of the Fresnel surfaces of
its characteristics in both convex and non-convex situations.
In Section
\ref{SEC:examples} we illustrate our assumptions
for differential hyperbolic systems, for hyperbolic
equations of the second order, and for general scalar
hyperbolic equations of higher orders.

\numberwithin{equation}{section}

\section{Prerequisites}
The first step is based on a decomposition of the extented phase
space $\R_+\times\R^n$ into zones, a pseudo-differential zone containing bounded frequencies and
a hyperbolic zone containing the large frequencies. It follows basically the treatment of \cite{RS05a}, 
\cite{Wir06}, going back to the papers of Reissig-Yagdjian for the treatment of wave equations with time-dependent propagation speed (\cite{RY00}, \cite{RY98} and references cited therein, and in \cite{RS05a}) and the book of 
Yagdjian, \cite{Yagdjian}, where the method of zones was introduced and applied to the investigation of the
Cauchy problem for hyperbolic operators with multiple characteristics.
 
\subsection{Hyperbolic symbol classes}
We make use of the implicitly defined function $t_\xi$ from
\begin{equation}\label{eq:zones:high-reg}
  (1+t_\xi)|\xi|=N (\log(e+t_\xi))^\nu
\end{equation}
with suitable constants $N$  and $\nu$ (usually required to belong to $[0,1]$) and define the zones
\begin{equation}\label{eq:zones:high-reg2}
  Z_{hyp}(N,\nu)=\{(t,\xi)|t\geq t_\xi\},\qquad\qquad Z_{pd}(N)=\{(t,\xi)|0\leq t\leq t_\xi\}.
\end{equation} 
In $Z_{hyp}(N,\nu)$ we apply a diagonalisation procedure to the full symbol. The basic idea of this 
diagonalisation scheme comes from the treatment of degenerate hyperbolic problems and is closely 
related to the approach of \cite{RS05a}. If $\nu>0$, the hyperbolic zone is subdivided further
into two parts, the regular and the oscillating sub-zones
\begin{equation}\label{eq:zones:high-reg3}
  Z_{reg}(N,\nu)=Z_{hyp}(N,2\nu) \qquad\qquad Z_{osc}(N)=Z_{hyp}(N,\nu)\setminus Z_{hyp}(N,2\nu).
\end{equation} 
We denote the new boundary as $\tilde t_\xi$. It satisfies  $(1+\tilde t_\xi)|\xi|=N\big(\log(e+\tilde t_\xi)\big)^{2\nu}$. The parameter $\nu$ is used to control the allowed amount of oscillations in the time-dependence of coefficients and amplitude functions. Based on the terminology from \cite[Section 3.11.1]{Yagdjian} and of \cite{RY98}, \cite{RY00a}, one speaks
of {\em very slow oscillations} if $\nu=0$, {\em slow oscillations} of $\nu\in(0,1)$, and {\em fast oscillations} if $\nu=1$. The following definition of symbol classes and also the basic assumptions introduced later on will depend on the parameter $\nu$.

\begin{defn}\label{def:Symb_hyp}
  The time-dependent Fourier multiplier $a(t,\xi)$ belongs to the {symbol class}
  $\Sc_{N,\nu}^{\ell_1,\ell_2}\{m_1,m_2\}$ if it satisfies the symbolic estimates
  \begin{equation}\label{eq:symb_est_hyp}
    \left|\D_t^k\D_\xi^\alpha a(t,\xi)\right|\leq 
    C_{k,\alpha} |\xi|_t^{m_1-|\alpha|}\left(\frac1{1+t} \big(\log(e+t)\big)^\nu \right)^{m_2+k}
  \end{equation}
  for all multi-indices 
  $\alpha\in\N^n$ with $|\alpha|\leq\ell_1$ and all natural numbers $k\leq\ell_2$, where 
  \begin{equation}
    |\xi|_t = \max\left(|\xi|,\frac N{1+t}\big(\log(e+t)\big)^\nu\right).
  \end{equation}
\end{defn}

If the symbolic estimates hold for all derivatives we write 
$\Sc_{N,\nu}\{m_1,m_2\}$ short for 
$\Sc_{N,\nu}^{\infty,\infty}\{m_1,m_2\}$. 
Furthermore, the definition extents immediately 
to matrix-valued 
Fourier multipliers. The rules of the corresponding 
symbolic calculus are simple consequences of 
Definition~\ref{def:Symb_hyp}
together with 
\eqref{eq:zones:high-reg}--\eqref{eq:zones:high-reg3} 
and are collected in the following proposition.

\begin{prop}\label{prop:II:calc_rules}
\begin{enumerate}
\item $\Sc_{N,\nu}^{\ell_1,\ell_2}\{m_1,m_2\}$ is a vector space.
\item $ \Sc_{N,\nu}^{\ell_1',\ell_2'}\{m_1-k,m_2+\ell\}\hookrightarrow \Sc_{N,\nu}^{\ell_1,\ell_2}\{m_1,m_2\}$ for all  
  $\ell\geq k\geq0$, $\ell_1'\geq\ell_1$ and $\ell_2'\geq \ell_2$.
\item $\Sc_{N,\nu}^{\ell_1,\ell_2}\{m_1,m_2\}\cdot \Sc_{N,\nu}^{\ell_1,\ell_2}\{m_1',m_2'\}\hookrightarrow \Sc_{N,\nu}^{\ell_1,\ell_2}\{m_1+m_1',m_2+m_2'\}$.
\item $\D_t^k\D_\xi^\alpha \Sc_{N,\nu}^{\ell_1,\ell_2}\{m_1,m_2\}\hookrightarrow \Sc_{N,\nu}^{\ell_1-|\alpha|,\ell_2-k}\{m_1-|\alpha|,m_2+k\}$.
\item $\Sc_{N,\nu}^{0,0}\{-1,2\}\hookrightarrow L^\infty_\xi L^1_t(Z_{reg}(N,\nu))$.
\end{enumerate}
\end{prop}
\begin{pf}
The statements are simple consequences of the symbolic estimates. Since it is vitally important for us later on, we show the last one. It follows that 
\begin{align*}
\int_t^\infty |a(\tau,\xi)| \d\tau &\lesssim \int_t^\infty \frac{(\log(e+\tau))^{2\nu}}{|\xi| (e+\tau)^2}\d\tau 
\\&= \left[  -\frac{(\log(e+\tau))^{2\nu}}{|\xi| (e+\tau)} \right]_t^\infty + 2\nu \int_t^\infty  \frac{(\log(e+\tau))^{2\nu-1}}{|\xi| (e+\tau)^2}\d\tau \\&\lesssim  \frac{(\log(e+t))^{2\nu}}{|\xi| (1+t)} \lesssim 1
\end{align*}
for any $a\in\Sc_{N,\nu}^{0,0}\{-1,2\}$ 
uniform in $(t,\xi)\in Z_{reg}(N,\nu).$ 
\end{pf}

Of particular importance are the embedding relations (2). 
They constitute a symbolic hierarchy, which is used in the diagonalisation scheme, cf.~Sect.~\ref{diag-scheme}. We define the residual symbol classes
\begin{equation}
  \H_{N,\nu}\{m\}=\bigcap_{k\in\Z} \Sc_{N,\nu}\{m-k,k\}
\end{equation}
as intersections along the embeddings 
$\Sc_{N,\nu}\{m_1-k,m_2+k\}\hookrightarrow \Sc_{N,\nu}\{m_1,m_2\}$ for $k\ge0$.

For later use we define the cut-off function $\chi_{pd}(t,\xi)=\chi(|\xi|(1+t)(\log(e+t))^{-\nu}/N)$ for $\chi\in C_0^\infty(\R)$, $\chi(s)=1/2$ for $s\le 1$ and $\chi(s)=0$ for $s>1$ and similarly  $\chi_{hyp}(t,\xi)=1-\chi_{pd}(t,\xi)$.

\subsection{Basic assumptions}
We collect our assumptions on the symbol $A(t,\xi)$. Throughout this 
paper we require

\medskip
{\bf (A1)$_{\ell_1,\ell_2}$} {\sl Hyperbolic operator of first order with bounded coefficients.}
We assume that the matrix operator $A(t,\D)$ has a symbol satisfying
\begin{equation}
A(t,\xi) \in \Sc_{N,\nu}^{\ell_1,\ell_2}\{1,0\}
\end{equation}
for some $N>0$ and $\nu\in[0,1]$, which has all its eigenvalues in a strip $\spec A(t,\xi)\subseteq S_c= \{\zeta\in\C:|\Im\zeta|\le c\}$ for all $(t,\xi)\in Z_{hyp}(N,\nu)$.  

We say the system is {\sl symmetric hyperbolic} 
if  we assume in addition
\begin{equation}
  \Im A(t,\xi) = \frac1{2\i}(A(t,\xi)-A^*(t,\xi)) \in  \Sc_{N,\nu}^{\ell_1,\ell_2}\{0,1\}.
\end{equation}
We write {\bf (A1$^+$)} for this stronger assumption.

\medskip
{\bf (A2)} {\sl Uniform strict hyperbolicity up to $t=\infty$.} 
We assume that there exists a homogeneous of order one in
$\xi$ matrix $A_1(t,\xi)$ with 
$A(t,\xi)-A_1(t,\xi)\chi_{hyp}(t,\xi)\in 
\Sc_{N,\nu}^{\ell_1,\ell_2}\{0,1\}$,
which has real and distinct eigenvalues 
$$\lambda_1(t,\xi),... , \lambda_m(t,\xi)\in\Sc_{N,\nu}^{\ell_1,\ell_2}\{1,0\},$$
the so-called {\em characteristic roots} of the symbol 
$A_1(t,\xi)$. In ascending order we 
denote them as $\lambda_1(t,\xi),... , \lambda_m(t,\xi)$. 
Furthermore, we assume that 
\begin{equation}
   \liminf_{t\to\infty} \min_{\omega\in\S^{n-1}} 
   |\lambda_i(t,\omega)-\lambda_j(t,\omega)| >0
\end{equation}
for all $i\ne j$.

As a consequence of the first two assumptions we get 
\begin{prop}\label{PROP:characteristics}
Assume (A1)$_{\ell_1,\ell_2}$ and (A2).
For all $j=1, ...,  m$ the characteristic roots satisfy $\lambda_j(t,\xi)\in\Sc_{N,\nu}^{\ell_1,\ell_2}\{1,0\}$ and
for all $i\ne j$ their difference satisfies $(\lambda_i(t,\xi)-\lambda_j(t,\xi))^{-1}\in\Sc_{N,\nu}^{\ell_1,\ell_2}\{-1,0\}$.
Furthermore, the eigenprojection $P_j(t,\xi)$ corresponding to $\lambda_j(t,\xi)$ satisfies
$P_j(t,\xi)\in\Sc_{N,\nu}^{\ell_1,\ell_2}\{0,0\}$.
\end{prop}
\begin{pf}[Sketch of proof.]
The properties of the characteristic roots follow from the spectral estimate $|\lambda_j(t,\omega)|\le||A_1(t,\omega)||$ together with
the obvious symbol properties of the coefficients of the characteristic polynomial and the uniform strict
hyperbolicity. The eigenprojections can be expressed in terms of the characteristic roots 
\begin{equation}
  P_j(t,\xi)=\prod_{i\ne j} \frac{A_1(t,\xi)-\lambda_i(t,\xi)\mathrm I}{\lambda_j(t,\xi)-\lambda_i(t,\xi)}
\end{equation}
and again the symbolic calculus yields the desired result.
\end{pf}

For $a_1,\ldots,a_m\in\C$ we define $\diag(a_1,\ldots,a_m)$
to be the diagonal matrix with entries $a_1,\ldots,a_m$ along the
diagonal. For a matrix $M$, we define $\diag M$ as the
diagonal matrix with $M_{11},\ldots,M_{mm}$ at the diagonal.

\begin{prop}\label{prop4} Assume (A1)$_{\ell_1,\ell_2}$ and (A2).
   There exists a uniformly bounded and invertible matrix 
   $M(t,\xi)\in\Sc_{N,\nu}^{\ell_1,\ell_2}\{0,0\}$ with 
   $M^{-1}(t,\xi)\in\Sc_{N,\nu}^{\ell_1,\ell_2}\{0,0\}$  
   which diagonalises the symbol $A_1(t,\xi)$,
   \begin{multline}
     A_1(t,\xi)M(t,\omega)=M(t,\xi)\mathcal D(t,\xi),\quad
     \mathcal D(t,\xi)=\diag\big(\lambda_1(t,\xi),...,\lambda_m(t,\xi)\big).
   \end{multline}
\end{prop}
\begin{pf}
By Proposition~\ref{PROP:characteristics} we see that the symmetriser of $A_1(t,\xi)$ given as
$H(t,\xi)=\sum_{j} P_j^*(t,\xi)P_j(t,\xi)$ satisfies $H(t,\xi)\in \Sc_{N,\nu}^{\ell_1,\ell_2}
\{0,0\}$. Our aim is to express $M^{-1}(t,\xi)$ explicitly in terms of  the eigenprojections 
$P_j(t,\xi)\in \Sc_{N,\nu}^{\ell_1,\ell_2}\{0,0\}$. Let therefore be $v_j(t,\xi)$ be (a smoothly chosen)
unit vector from the one-dimensional $j$-th eigenspace $\mathrm{ran} \,P_j(t,\xi)$ and  
$M^{-1}(t,\xi) u $ be the vector with the entries
given by the scalar products $(v_j, P_j(t,\xi) u)=(v_j, H(t,\xi) u)$  for any vector 
$u\in\mathbb C^n$. 

Note that $v_j(t,\xi)$ is unique up to a sign and locally in $t$ and $\xi$ expressable as $v_j(t,\xi) = P_j(t,\xi) e / \|P_j(t,\xi) e\|$ for some fixed unit vector $e$ chosen away from the orthogonal complement to the eigenspace. Differentiating this yields the desired bounds 
$v_j(t,\xi)\in\Sc_{N,\nu}^{\ell_1,\ell_2} \{ 0,0\}$ and therefore $M^{-1}(t,\xi)\in\Sc_{N,\nu}^{\ell_1,\ell_2}\{0,0\}$. The inverse matrix $M(t,\xi)$ has the vectors $v_j(t,\xi)$ as its columns
and therefore $M(t,\xi)\in \Sc_{N,\nu}^{\ell_1,\ell_2} \{0,0\}$ as claimed.
\end{pf}

\begin{rem}
Under Assumption (A1$^+$) the matrix $M(t,\xi)$ can 
be chosen unitary, which simplifies some of the 
considerations later on.
\end{rem}

\medskip
{\bf (A3)} {\sl Generalised energy conservation property.} 
We have to make a technical assumption which guarantees 
us that the lower order terms are not too strong. 
This can be written as
\begin{equation}\label{eq:1.9}
   \sup_{(s,\xi),(t,\xi)\in Z_{hyp}(N,\nu)}  \left\| \int_s^t 
   \Im F_0(\theta,\xi) \d\theta \right\|  < \infty
\end{equation}
for the matrix 
\begin{equation}
F_0(t,\xi) = \diag( M^{-1}(t,\xi) (A(t,\xi)-A_1(t,\xi)) 
M(t,\xi) - M^{-1}(t,\xi)(\D_t M(t,\xi))).
\end{equation}

\begin{rem}
This statement does not follow from symbol 
estimates as we only have 
$F_0\in \Sc_{N,\nu}^{\ell_1,\ell_2-1}\{0,1\}$. We will later use it 
in Theorem \ref{THM:enegry} relating it to 
uniform bounds on the micro-energy of solutions within 
$Z_{reg}(N,\nu)$ explaining the nature of (A3). 
Assumption (A3) is independent 
of the choice of the diagonaliser $M(t,\xi)$, provided it satisfies the properties of Proposition~\ref{prop4}.
\end{rem}

If we assume (A1$^+$), then the diagonaliser $M(t,\xi)$ can be chosen unitary. This implies
\begin{prop}
Assume (A1$^+$).  Then the matrix $M^*(\D_t M)$ 
is self-adjoint and therefore
\begin{equation}
      \Im\ \diag M^*(t,\xi)\D_t M(t,\xi)=0,
\end{equation}
in particular, $\Im\ F_0(t,\xi) = \Im\
\diag( M^*(t,\xi) (A(t,\xi)-A_1(t,\xi)) M(t,\xi)).$
\end{prop}
\begin{pf}
Obviously $0=\D_t \mathrm I=(\D_t M^*)M+M^*(\D_t M)$ and therefore 
$$(M^*\D_t M)^*=(\D_t M)^* M=-(\D_t M^*) M=M^*(\D_tM).$$ But this implies that the diagonal entries of $(\D_tM^*)M$ are real.
\end{pf}

Assumption (A3) can be replaced by the weaker assumption that there exists a real scalar function $\delta(t,\xi)\in \Sc\{0,1
\}$ such that $\Im F_0(t,\xi) - \delta(t,\xi) I$ satisfies the equivalent of \eqref{eq:1.9}, i.e.
\begin{equation}
   \sup_{(s,\xi),(t,\xi)\in Z_{reg}(N,\nu)}  \left\| \int_s^t (\Im F_0(\theta,\xi) - \delta(\theta,\xi) I) \d\theta \right\|  < \infty.
\end{equation}
We refer to this weaker assumption as {\bf (A3$^-$)}.

\medskip
{\bf (A4)} {\sl Weak dissipativity.} There exists a constant $c$ and a function $\gamma(t)$ such that
\begin{equation}
  \Im A(t,\xi) + c|\xi| \mathrm I  + \textstyle \frac12 \gamma(t)\mathrm I \ge 0
\end{equation}
in the sense of matrices and within $Z_{pd}(N,\nu)$ for sufficiently large $N$, where
\begin{equation}\label{eq:hCond}
  \int_0^t \gamma(s) \d s \le C (\log(e+t))^\nu
\end{equation}
for some constant $C$ and all $t\ge 0$.

If we assume (A3$^-$), we also weaken (A4) to {\bf (A4$^-$)}
\begin{equation}
      \Im A(t,\xi) + c|\xi|\mathrm I  + \textstyle \frac12 \gamma(t) \mathrm I + \frac12 \delta(t,\xi) \mathrm I \ge 0
\end{equation}
with the same $\delta(t,\xi)$ as in (A3$^-$).

\begin{rem}
This assumption is mainly for convenience and to simplify the main estimate in the pseudo-differential zone given in Lemma~\ref{lemPD1}. If one wants to obtain sharp estimates in specific situations a more detailed analysis of the behaviour of solutions for small frequencies has to be carried out. We refer to \cite{Wir06} for an example, where such a more detailed consideration was essential.  

{ 
Another suitable replacement for (A4) could be to assume that the full symbol $A(t,\xi)$ has all its eigenvalues in the (closed) complex upper half plane and that it is uniformly symmetrisable on $Z_{pd}(N',\nu)$. 
}
\end{rem}

Later on we will discuss the construction of representations of solutions and resulting estimates under assumptions 
(A1) to (A4). Under the weaker assumptions (A3$^-$) and (A4$^-$) the uniform behaviour of the energy within $Z_{reg}(N,\nu)$ is replaced by an equivalence to $\exp\left(\int_0^t \delta(\tau,\xi)\d\tau\right)$ and can be obtained by minor modifications in the arguments.

\section{Representation of solutions}
In order to solve \eqref{eq:CP} we apply a partial Fourier transform with respect to the spatial
variables. Assuming $U_0\in\mathcal S(\R^n)$, this gives a system of ordinary differential equations
\begin{equation}
   \D_t \widehat U = A(t,\xi)\widehat U,\qquad \widehat U(0,\cdot)=\widehat U_0
\end{equation}
parametrised by the frequency parameter $\xi\in\R^n$ with initial data $\widehat U_0\in\mathcal S(\R^n)$. We construct its fundamental solution
$\mathcal E(t,s,\xi)$, i.e. the solution to the matrix-valued problem
\begin{equation}\label{eq:Eeq}
  \D_t\mathcal E(t,s,\xi)=A(t,\xi)\mathcal E(t,s,\xi),\qquad \mathcal E(s,s,\xi)=\mathrm I\in\C^{m\times m}.
\end{equation}
Then, the inverse partial Fourier transform gives $U(t,x)=\mathcal E(t,0,\D)U_0(x)$ for
all $t$ and $x$, while our construction gives a representation of $\mathcal E(t,0,\D)$ as a matrix of
Fourier integral operators.

\subsection{Treatment in the pseudo-differential zone}\label{sec2.1}
We use the dissipativity assumption (A4) and restrict considerations to the zone $Z_{pd}(N,\nu)$. 
\begin{lem}\label{lemPD1}
Assume (A4). Then the fundamental solution to \eqref{eq:Eeq} satisfies
\begin{equation}
   \|\mathcal E(t,0,\xi)\| \le C \exp\left(C'\big(\log(e+t_\xi)\big)^\nu\right)
\end{equation}
for some constants $C$ and $C'$ and uniform in $(t,\xi)\in Z_{pd}(N,\nu)$.
\end{lem}
\begin{pf}
Fix $\xi\in\R^n$ and denote by $V(t)$ the solution to $\D_t V(t) = A(t,\xi) V(t)$ with data $V(0)=V_0$. Then using the Euclidean inner product $(\cdot,\cdot)$ on $\C^m$ we obtain from (A4) that
\begin{align}
\frac{\d}{\d t} \|V(t)\|^2 = -2\big(\Im A(t,\xi) V(t), V(t) \big) \le 2 c |\xi| \|V(t)\|^2 + \gamma(t) \|V(t)\|^2
\end{align}
for all $t\le t_\xi$. Hence, by applying the Gronwall inequality the estimate
\begin{equation}
  \|V(t)\|^2 \lesssim \|V_0\|^2 \exp\left(2ct|\xi| + \int_0^t \gamma(s)\d s\right) \lesssim \exp\left(C'\big(\log(e+t_\xi)\big)^\nu\right)
\end{equation}
follows. This proves the statement.
\end{pf}

Later on we need estimates for derivatives  of $\mathcal E(t_\xi,0,\xi)$ with respect to $\xi$ as $|\xi|\to0$. They essentially follow from the symbol estimates satisfied by $A(t,\xi)$ due to (A1)$_{\ell_1,\ell_2}$.
\begin{lem} \label{lemPD2}
Assume (A1)$_{\ell_1,\ell_2}$ and (A4).
For all multi-indices $|\alpha|\le \min(\ell_1,\ell_2)$ the estimate
\begin{equation}
  || \D_\xi^\alpha \mathcal E(t_\xi,0,\xi) || \le C_\alpha |\xi|^{-|\alpha|} \big(\log(e+t_\xi)\big)^{\nu|\alpha|} \exp\left(C'\big(\log(e+t_\xi)\big)^\nu\right)
\end{equation}
is valid uniformly on $|\xi|\le N$.
\end{lem}
\begin{pf} 
To prove this fact we use Duhamel's formula for 
$\xi$-derivatives of \eqref{eq:Eeq}. Let first $|\alpha|=1$.
Then $\D_t\D_\xi^\alpha\mathcal E=(\D_\xi^\alpha A)\mathcal E+
A(\D_\xi^\alpha\mathcal E)$ and thus
using $\D_\xi^\alpha\mathcal E(0,0,\xi)=0$ we obtain the representation
\begin{equation}
  \D_\xi^\alpha \mathcal E(t,0,\xi)=\int_0^t 
  \mathcal E(t,\tau,\xi)(\D_\xi^\alpha A(\tau,\xi))
  \mathcal E(\tau,0,\xi)\d\tau.
\end{equation}
Since $||\D_\xi^\alpha A(\tau,\xi)||\lesssim 1$, this
implies the estimate
\begin{align}
||\D_\xi^\alpha\mathcal E(t,0,\xi)||&\lesssim t\exp\left(C\big(\log(e+t)\big)^\nu\right)\\
&\lesssim |\xi|^{-1}\big(\log(e+t_\xi)\big)^\nu\exp\left(C\big(\log(e+t_\xi)\big)^\nu\right)\notag
\end{align} uniformly on $Z_{pd}(N,\nu)$.

For $|\alpha|=\ell>1$ we use Leibniz formula to represent 
$\D_\xi^\alpha \mathcal E(t,0,\xi)$ by a corresponding 
Duhamel integral using a sum of terms 
$(\D_\xi^{\alpha_1} A(t,\xi))(\D_\xi^{\alpha_2}
\mathcal E(t,0,\xi))$ for $|\alpha_1|+ |\alpha_2|\le\ell$, 
$\alpha_2<\alpha$, and apply induction over $\ell$ to obtain
\begin{equation}
   \| \D_\xi^\alpha \mathcal E(t,0,\xi) \| 
   \lesssim |\xi|^{-|\alpha|} 
   \big(\log(e+t_\xi)\big)^{\nu|\alpha|} 
   \exp\left(C'\big(\log(e+t_\xi)\big)^\nu\right),
   \qquad |\xi|\le N.
\end{equation}

Derivatives of $\mathcal E$ with respect to $t$ are 
estimated directly by the differential equation and 
$\|A(t,\xi)\|\lesssim |\xi|$. Combined with  
\begin{equation}\label{eq:2.10}
|\D^\alpha_{\xi} t_\xi| \lesssim t_\xi |\xi|^{-\alpha} 
\lesssim |\xi|^{-1-|\alpha|} \big(\log(e+t_\xi)\big)^\nu
\end{equation}
this completes the proof.
\end{pf}

\begin{rem} If one is only interested in energy estimates, the statement of Lemma~\ref{lemPD1} is sufficient to get
the corresponding micro-energy estimate in the pseudo-differential zone. 
\end{rem}
\begin{rem}
Any replacement of (A4), which is sufficient to deduce Lemma~\ref{lemPD1}, is also enough to prove Lemma~\ref{lemPD2}.
This is interesting in particular for the treatment for second order scalar problems, where other
ways to prove the estimate of Lemma 9 are available and important to obtain sharp result. See, e.g., the treatment in
\cite{Wir06} or \cite{HW08} where reformulations as integral equations were used.
\end{rem}

\subsection{Treatment in the hyperbolic zone}
The following considerations are only applied within 
the hyperbolic zone. We can write formulae globally 
for all $t$ and $\xi$ if we make use of the cut-off 
function $\chi_{hyp}(t,\xi)$. We omit this to keep 
notation as simple as possible.

\subsubsection{Treatment of the main part}
We apply the diagonaliser $M(t,\xi)$ of $A_1(t,\xi)$. 
For this we set $U^{(0)}=M^{-1}(t,\xi)\widehat U$,
so that
\begin{align}
   D_t U^{(0)}
   &=(\D_tM^{-1}(t,\xi))\widehat U+M^{-1}(t,\xi)\D_t 
   \widehat U\notag\\&=\big((\D_tM^{-1})M+ M^{-1}(t,\xi)A(t,\xi)M(t,\xi)\big)U^{(0)}\\
   &=\big(\mathcal D(t,\xi)+R_0(t,\xi) \big)U^{(0)}\notag
\end{align}
is valid. Now the main term $\mathcal D(t,\xi)=
\diag(\lambda_1(t,\xi),\ldots ,\lambda_m(t,\xi))$ is diagonal and
the remainder is given by
\begin{equation}
  R_0(t,\xi) = M^{-1}(t,\xi) (A(t,\xi)-A_1(t,\xi)) M(t,\xi) + (\D_tM^{-1}(t,\xi))M(t,\xi).
\end{equation}


The rules of the symbolic calculus imply  
$R_0(t,\xi)\in\Sc_{N,\nu}^{\ell_1,\ell_2-1}\{0,1\}$, 
while our assumptions on the characteristic roots give 
$\mathcal D(t,\xi)\in\Sc_{N,\nu}^{\infty,\ell_2}\{1,0\}$.  

\subsubsection{The diagonalisation scheme}\label{diag-scheme} 
 The initial diagonalisation step was done in the 
 previous section. We transformed the system
 by a matrix $M\in\Sc_{N,\nu}^{\infty,\ell_2}\{0,0\}$ to a diagonal form modulo $\Sc_{N,\nu}^{\ell_1,\ell_2-1}\{0,1\}$. In the $k$-th step
 we want to diagonalise it modulo the class
 $\Sc_{N,\nu}^{\ell_1,\ell_2-k-1}\{-k,k+1\}$ for sufficiently large $N$.
 (The choice of $N$  depends on $k$ in order to guarantee the invertibility of the transformation matrix $N_k\in\Sc_{N,\nu}\{0,0\}$ below.)
 
Recursively, we construct matrices $N_k(t,\xi)$ to diagonalise the system $\D_t-\mathcal D(t,\xi)-R_0(t,\xi)$ modulo 
$\Sc_{N,\nu}^{\ell_1,\ell_2-k-1}\{-k,k+1\}$. We use the notation
\begin{equation}
 N_k(t,\xi)= \sum_{j=0}^k N^{(j)}(t,\xi),\qquad F_k(t,\xi)=\sum_{j=0}^k F^{(j)}(t,\xi), 
\end{equation}
where initially $N^{(0)}(t,\xi)=\mathrm I$, $B^{(0)}(t,\xi)=R_0(t,\xi)$ and $F^{(0)}(t,\xi)=\mathrm{diag}\;B^{(0)}(t,\xi)$. The construction goes along the following recursive scheme. 

{\bf Step $k$, $k\geq1$.} We set
\begin{align}
      F^{(k-1)}(t,\xi)&=\mathrm{diag }\;B^{(k-1)}(t,\xi),\\
      \big(N^{(k)}(t,\xi)\big)_{i,j} &= \begin{cases} \frac{-B^{(k-1)}(t,\xi)_{ij}}{\lambda_i(t,\xi)-\lambda_j(t,\xi)}, \qquad&{i\ne j}\\ 0, & i=j, \end{cases} \\
      B^{(k)}(t,\xi) &= (\D_t-\mathcal D(t,\xi)-R_0(t,\xi))N_{k}(t,\xi)\notag\\&\qquad-N_{k}(t,\xi)(\D_t-\mathcal D(t,\xi)-F_{k-1}(t,\xi)).
\end{align}

We prove by induction over $k$ the estimates $N^{(k)}(t,\xi)\in \Sc_{N,\nu}^{\ell_1,\ell_2-k}\{-k,k\}$ and 
$B^{(k)}(t,\xi)\in \Sc_{N,\nu}^{\ell_1,\ell_2-k-1}\{-k,k+1\}$.  For $k=1$ we know
\begin{multline}
  F^{(0)}(t,\xi)\in \Sc_{N,\nu}^{\ell_1,\ell_2-1}\{0,1\}, \\
  N^{(1)}(t,\xi)\in \Sc_{N,\nu}^{\ell_1,\ell_2-1}\{-1,1\},\\
  B^{(1)}(t,\xi)\in \Sc_{N,\nu}^{\ell_1,\ell_2-2}\{-1,2\},
\end{multline}
the last one follows from the representation $$B^{(1)}(t,\xi)=\D_t N^{(1)}(t,\xi)-(R_0(t,\xi)-F^{(0)}(t,\xi))N^{(1)}(t,\xi).$$ 

For $k\geq1$, assume we know $B^{(k-1)}(t,\xi)\in \Sc_{N,\nu}^{\ell_1,\ell_2-k}\{-k+1,k\}$. 
Then by definition of $N^{(k)}$ we have from $(\lambda_i(t,\xi)-\lambda_j(t,\xi))^{-1}\in \Sc_{N,\nu}^{\infty,\ell_2}\{-1,0\}$ 
that $N^{(k)}(t,\xi)\in \Sc_{N,\nu}^{\ell_1,\ell_2-k}\{-k,k\}$ and obviously $F^{(k-1)}(t,\xi)\in \Sc_{N,\nu}^{\ell_1,\ell_2-k}\{-k+1,k\}$. 
Moreover, by algebraic manipulations we get
\begin{align*}
     B^{(k)}&=(\D_t-\mathcal D-R_0)(\sum_{\nu=0}^{k} N^{(\nu)})
        -(\sum_{\nu=0}^{k} N^{(\nu)})(\D_t-\mathcal D-\sum_{\nu=0}^{k-1} F^{(\nu)})\\
     &=B^{(k-1)}+[N^{(k)},\mathcal D]-F^{(k-1)}+\D_tN^{(k)}+R_0N^{(k)}\\
      &\qquad\qquad\qquad\qquad\qquad\qquad+N^{(k)}\sum_{\nu=0}^{k-1} F^{(\nu)} 
    -(\sum_{\nu=1}^{k}N^{(\nu)})F^{(k)}.
\end{align*}
The definition of $N^{(k)}(t,\xi)$ was chosen in such a way that we have 
\begin{equation}
 B^{(k-1)}(t,\xi)+[N^{(k)}(t,\xi),\mathcal D(t,\xi)]-F^{(k-1)}(t,\xi)=0.
\end{equation}
The sum of the remaining terms belongs to the symbol class $\Sc_{N,\nu}^{\ell_1,\ell_2-k-1}\{-k,k+1\}$. Hence
we have proved $B^{(k)}(t,\xi)\in \Sc_{N,\nu}^{\ell_1,\ell_2-k-1}\{-k,k+1\}$.
Furthermore, the definition of $B^{(k)}(t,\xi)$ implies the operator identity
\begin{equation}\label{eq:II:perf_diag}
  \big(\D_t-\mathcal D(t,\xi)-R_0(t,\xi)\big)N_k(t,\xi)
  =N_k(t,\xi)\big(\D_t-\mathcal D(t,\xi)-F_{k-1}(t,\xi)\big)
\end{equation}
modulo $\Sc_{N,\nu}^{\ell_1,\ell_2-k-1}\{-k,k+1\}.$
Thus we have constructed the desired diagonaliser if we can show that the matrix
$N_k(t,\xi)$ is invertible on $Z_{hyp}(N,\nu)$ with uniformly bounded inverse. But
this follows from $N_k-I\in \Sc_{N,\nu}^{\ell_1,\ell_2-k}\{-1,1\}$ by the choice of a sufficiently large 
zone constant $N$. Indeed, we have
\begin{equation}
 ||N_k(t,\xi)-I|| \leq C \frac{(\log(e+t))^\nu}{|\xi|(1+t)}\leq \frac{C'}{N}\to0
   \qquad\text{as $N\to\infty$}. 
\end{equation}
Thus with the notation $R_k(t,\xi)=-N_k^{-1}(t,\xi)B^{(k)}(t,\xi)$ we have proven the 
following lemma.

\begin{lem}\label{lem:II:perf_diag}
 Assume (A1)$_{\ell_1,\ell_2}$ and (A2).  For each $1\leq k< \ell_2$ there exists a zone constant $N$ and matrix valued symbols
  \begin{itemize}
  \item $N_k(t,\xi)\in \Sc_{N,\nu}^{\ell_1,\ell_2-k}\{0,0\}$, which is invertible for all $(t,\xi)\in Z_{hyp}(N,\nu)$ 
    and its inverse satisfies $N_k^{-1}(t,\xi)\in \Sc_{N,\nu}^{\ell_1,\ell_2-k}\{0,0\}$;
  \item $F_{k-1}(t,\xi)\in \Sc_{N,\nu}^{\ell_1,\ell_2-k}\{0,1\}$, which is diagonal;
  \item $R_k(t,\xi)\in \Sc_{N,\nu}^{\ell_1,\ell_2-k-1}\{-k,k+1\}$; 
  \end{itemize}
  such that the (operator) identity
  \begin{multline}
  \big(\D_t-\mathcal D(t,\xi)-R_0(t,\xi)\big)N_k(t,\xi)\\
  =N_k(t,\xi)\big(\D_t-\mathcal D(t,\xi)-F_{k-1}(t,\xi)-R_k(t,\xi)\big)
  \end{multline}
  holds true for all $(t,\xi)\in Z_{hyp}(N,\nu)$.
\end{lem}

\subsubsection{Remarks on perfect diagonalisation} 
Under assumption (A1)$_{\infty,\infty}$ Lemma \ref{lem:II:perf_diag} can be 
understood as a perfect diagonalisation of the original system. 
If we define $F(t,\xi)$ as an asymptotic sum of the $F^{(k)}(t,\xi)$, 
\begin{equation}
  F(t,\xi)\sim \sum_{k=0}^\infty F^{(k)}(t,\xi),
\end{equation}
this means we require $F(t,\xi)-F_k(t,\xi)\in \Sc_{N,\nu}\{-k-1,k+2\}$ for all $k\in\N$,
and similarly
\begin{equation}
  N(t,\xi)\sim\sum_{k=0}^\infty N^{(k)}(t,\xi),
\end{equation}
which can be chosen to be invertible,
equation \eqref{eq:II:perf_diag} implies
\begin{equation}
  \big(\D_t-\mathcal D(t,\xi)-R(t,\xi)\big)N(t,\xi)-N(t,\xi)\big(\D_t-\mathcal D(t,\xi)-F(t,\xi)\big)
  \in \H_{N,\nu}\{1\}.
\end{equation}
In this sense the symbol class $\H_{N,\nu}\{1\}$ is understood as the residual symbol class within the 
scheme of diagonalisation. 

\subsubsection{Solving the diagonalised system in the regular sub-zone}
Our strategy is as follows. In the first step we solve the diagonal system $\D_t-\mathcal D-F_{k-1}$. Its
fundamental solution $\widetilde{\mathcal E_k}(t,s,\xi)$ is given by
\begin{align}\label{eq:Ekdef}
\widetilde{\mathcal E_k}(t,s,\xi)
&=\exp\left(\i\int_s^t \big(\mathcal D(\tau,\xi)+F_{k-1}(\tau,\xi)\big)\d\tau\right)\notag\\
&=\exp\left(\i\int_s^t \mathcal D(\tau,\xi)\d\tau\right)\; \exp\left(\i\int_s^t F_{k-1}(\tau,\xi)\d\tau \right).
\end{align}
The second factor is uniformly bounded with uniformly bounded inverse on $Z_{reg}(N,\nu)$, because $F_{k-1}-F^{(0)}\in\Sc_{N,\nu}^{0,0}\{-1,2\}\subseteq L^\infty_\xi L^1_t(Z_{reg}(N,\nu))$ by Proposition
\ref{prop:II:calc_rules}, (5), and $\Im F^{(0)}$ satisfies (A3). 
The first factor is a unitary matrix. Thus, we obtain the uniform bound
\begin{equation}\label{eq:2.26}
\left\|  \ \widetilde{\mathcal E_k}(t,s,\xi) \right\|
 \le C,\qquad (t,\xi),(s,\xi)\in Z_{reg}(N,\nu)
\end{equation}
regardless of the order of $s$ and $t$.

In the second step we make the {\it ansatz}  
\begin{equation}
   \mathcal E_k(t,s,\xi)=\widetilde{\mathcal E_k}(t,s,\xi)\mathcal Q_k(t,s,\xi)
\end{equation}
for the fundamental solution $\mathcal E_k(t,s,\xi)$ to $\D_t-\mathcal D-F_k-R_k$. Then the matrix
$\mathcal Q_k(t,s,\xi)$ satisfies
\begin{equation}\label{eq:Qeq}
  \D_t \mathcal Q_k(t,s,\xi)=\mathcal R_k(t,s,\xi)\mathcal Q_k(t,s,\xi),\qquad \mathcal Q_k(s,s,\xi)=I\in \C^{m\times m}
\end{equation}
with coefficient matrix
\begin{equation}
  \mathcal R_k(t,s,\xi)=\widetilde{\mathcal E_k}(s,t,\xi)R_k(t,\xi)\widetilde{\mathcal E_k}(t,s,\xi).
\end{equation}
$R_k\in\Sc_{N,\nu}\{-k,k+1\}$ with $k\ge1$ implies uniform integrability of $\mathcal R_k$ over the regular sub-zone. Therefore, the representation of $\mathcal Q_k(t,s,\xi)$ as the Peano-Baker series 
\begin{multline}\label{eq:Qkrep}
  \mathcal Q_k(t,s,\xi)=\mathrm I+\sum_{k=1}^\infty \i^k \int_s^t \mathcal R_k(t_1,s,\xi)\int_s^{t_1}
  \mathcal R_k(t_2,s,\xi)\\\cdots\int_s^{t_{k-1}}\mathcal R_k(t_k,s,\xi)\d t_k\cdots\d t_1
\end{multline}
implies boundedness of $\mathcal Q_k(t,s,\xi)$ uniformly in $(t,\xi),(s,\xi)\in Z_{reg}(N,\nu)$.
We obtain even a little bit more:

\begin{lem}\label{lem5}
Assume (A1)$_{0,\ell_2}$, (A2) and (A3). Then the matrices $\mathcal Q_k(t,s,\xi)$, $k < \ell_2$, defined as solutions to \eqref{eq:Qeq}
are uniformly bounded on $Z_{hyp}(N,\nu)$, invertible and tend locally uniformly in $\xi\ne0$ to $\mathcal Q_k(\infty,s,\xi)\in L^\infty(Z_{hyp}(N,\nu))$ as $t\to\infty$.
\end{lem}
\begin{pf}
The convergence follows directly by applying the Cauchy criterion to the representation \eqref{eq:Qkrep}.
Invertibility is a consequence of the Liouville theorem, which gives
\begin{align}
\det\mathcal Q_k(t,s,\xi)&=\exp\left(\mathrm i\int_s^t \mathrm{trace}\,\mathcal R_k(\theta,s,\xi)\d\theta \right)= \exp\left(\mathrm i \int_s^t \mathrm{trace}\, R_k(\theta,\xi)\d\theta\right)\notag\\
&
\ge \exp\left(-m\int_s^t \|R_k(\theta,\xi)\|\d\theta\right)\ge C > 0,
\end{align}
uniform in $(t,\xi),(s,\xi)\in Z_{hyp}(N,\nu)$.
\end{pf}

The following theorem explains the true nature of 
assumption (A3) provided that $\nu=0$. 

\begin{thm} \label{THM:enegry}
Assume (A1)$_{0,\infty}$ with $\nu=0$ and (A2). 
Then (A3) is equivalent to the existence of a constant $C>0$ such that
\begin{equation}
   C^{-1} \|V\| \le \| \mathcal E(t,s,\xi) V\| \le C \|V\|
\end{equation}
for all vectors $V\in \C^m$ and $(t,\xi),(s,\xi)\in Z_{hyp}(N,0)$.
\end{thm}
\begin{pf}
Lemma~\ref{lem5} in combination with \eqref{eq:2.26}  gives the uniform bound under (A3). 
Without (A3) the estimate of \eqref{eq:2.26} has to be replaced by a polynomial bound
\begin{equation}
  \|\widetilde{\mathcal E}_k(t,s,\xi)\|, \|\widetilde{\mathcal E}_k(s,t,\xi)\| \le C_k \left(\frac{1+t}{1+s}\right)^K ,\qquad t\ge s,
\end{equation}
where the constant $K$ is independent of $k$ (and depends only on the first diagonal lower order term $F^{0}$). Similarly, we obtain with the same $K$
\begin{equation}
  \|{\mathcal E}_k(t,s,\xi)\| \le \exp\left( \int_s^t \| \Im (F_{k-1}(\tau,\xi)+R_k(\tau,\xi))\| \d\tau\right)  \le C_k' \left(\frac{1+t}{1+s}\right)^K 
\end{equation}
for all $t\ge s$. Choosing $k$ large, the polynomial decay of the remainder $R_k(t,\xi)$ becomes strong enough to compensate the increasing terms and we obtain the Duhamel representation
\begin{equation}\label{eq:2.16}
   \mathcal E_k(t,s,\xi) = \widetilde{\mathcal E}_k(t,s,\xi) \mathcal Z_k(s,\xi) - \i \int_t^\infty  \widetilde{\mathcal E}_k(t,\theta,\xi) R_k(\theta,\xi) \mathcal E_k(\theta,s,\xi)\d \theta
\end{equation}
with 
\begin{equation}
  \mathcal Z_k(s,\xi) = \mathrm I + \i \int_s^\infty \widetilde{\mathcal E}_k(s,\theta,\xi) R_k(\theta,\xi) \mathcal E_k(\theta,s,\xi)\d \theta \lesssim 1.
\end{equation}
Indeed, the integral in \eqref{eq:2.16} is convergent for $k> 2K$,
\begin{multline*}
\left\|\int_t^\infty  \widetilde{\mathcal E}_k(t,\theta,\xi) R_k(\theta,\xi) \mathcal E_k(\theta,s,\xi)\d \theta \right\|
\\\lesssim |\xi|^{-k} \int_t^\infty \left(\frac{1+\theta}{1+t}\right)^K  \frac1{(1+\theta)^{k+1}} \left(\frac{1+\theta}{1+s}\right)^K
\d\theta 
\end{multline*}
and bounded by $|\xi|^{-k}(1+t)^{K-k}(1+s)^{-K}$.
Similarly $\|\mathcal Z_k(s,\xi)-\mathrm I\| \le |\xi|^{-k}(1+s)^{-k}$ 
and hence the first term has the lower norm bound $\|\widetilde{\mathcal E}_k(t,s,\xi) \|\ge (1+t)^K(1+s)^{-K}$ for fixed $|\xi|$. Choosing $s$ big enough depending on $|\xi|$ implies that $\mathcal E_k(t,s,\xi)$ is a small perturbation of $\widetilde{\mathcal E}_k(t,s,\xi)$.

Assume now that (A3) is violated. Then we find sequences $t_\mu\to\infty$, $s_\mu$, and $\xi_\mu$ such that one matrix entry of the integral in \eqref{eq:1.9} tends to either $\infty$ or $-\infty$.  We consider the $+\infty$ case,  and assume w.l.o.g.~that $s_\mu>s$ for sufficiently big $s$ and that the matrix entry corresponds to the first diagonal element. Then with $e_1$ the first basis vector  $ \widetilde{\mathcal E}_k(t_\mu,s_\mu,\xi_\mu) e_1 \to \infty$ and therefore also $\mathcal E(t_\mu,s_\mu, \xi_\mu)e_1\simeq \mathcal E_k(t_\mu,s_\mu,\xi_\mu) N_k(s_\mu,\xi_\mu) M(s_\mu,\xi_\mu) e_1 \to \infty$ which contradicts to the uniform upper bound. Similarly, the $-\infty$ case contradicts to the lower bound and the statement is proven.
\end{pf}

\begin{rem}
A similar argument does not apply if $\nu>0$. In this situation the polynomial bound becomes a superpolynomial bound 
by $\exp(C(\log t)^{1+\nu})$ which cannot be compensated within the diagonalisation hierarchy. 
\end{rem}

The following results explain why we performed $k>1$ steps of diagonalisation. In order to prove dispersive estimates later on, we have to control derivatives of $\mathcal Q_k(t,x,\xi)$. Formally differentiating \eqref{eq:Qkrep} includes $\xi$-derivatives of $\widetilde{\mathcal E}_k(t,s,\xi)$, which are increasing in $t$. This can be compensated by better estimates on the remainder $R_k(t,\xi)$. As we will see below, more steps of the diagonalisation hierarchy allow for symbol-like estimates for more derivatives of the entries of  $\mathcal Q_k(t,x,\xi)$. 
  
\begin{lem}\label{lem:Zreg}
Assume (A1)$_{\ell_1,\ell_2}$, (A2) and (A3). Then the matrix $\mathcal Q_k(t,s,\xi)$ with $2k\le \ell_2$, satisfies for all $|\alpha|\le\min( k-1,\ell_1)$ 
  \begin{equation}
      \left\|\D_\xi^\alpha \mathcal Q_k(t,s,\xi)\right\| \lesssim |\xi|^{-|\alpha|}  \big(\log(e+|\xi|^{-1})\big)^{|\alpha|}
  \end{equation}
  uniformly in $s,t\ge \tilde t_\xi$. Furthermore, for $|\xi|\le N$ and $|\alpha|\le\min( \frac{k-1}2,\ell_1)$, we have
  \begin{equation}
      \left\|\D_\xi^\alpha \mathcal Q_k(t,\tilde t_\xi,\xi)\right\| \lesssim |\xi|^{-|\alpha|} \big(\log(e+|\xi|^{-1})\big)^{|\alpha|}
  \end{equation}
  uniformly in $t\ge \tilde t_\xi$.
\end{lem}
\begin{pf} 
The statement follows similar to the considerations in \cite[Lemma 2.10]{RS05a}. We concentrate on the second estimate, the first one is analogous. Note first, that by differentiating 
\eqref{eq:Ekdef} we obtain 
\begin{equation}
 \| \D_\xi^\alpha \widetilde{\mathcal E}_k(t, \tilde t_\xi,\xi) \| \lesssim (1+t)^{|\alpha|}
 \big(\log(e+t)\big)^{|\alpha|}, 
\end{equation}
 the logarithmic term arises from differentiating $F_{k-1}(t,\xi)$ in \eqref{eq:Ekdef}. This implies that $\mathcal R_k(t,t_\xi,\xi)$ satisfies weaker estimates than $R_k(t,\xi)$, for which we have $\Sc_N^{\ell_1,\ell_2-k-1}\{-k,k+1\}$ estimates.
Based on the equivalent to \cite[Proposition 11]{Wir06} we obtain 
\begin{equation}
 \| \D_\xi^\alpha \mathcal R_k(t,t_\xi,\xi) \| \lesssim |\xi|^{-1-|\alpha|} (1+t)^{-2}\big(\log(e+t)\big)^{|\alpha|}
\end{equation}
for $|\alpha|\le \frac{k-1}2$ and $|\alpha|\le \ell_1$. Differentiating the series representation \eqref{eq:Qkrep} term by term, taking into account these estimates and combining them with 
\begin{equation}\label{eq:est-txi-tilde}
|\D^\alpha_{\xi} \tilde t_\xi| \lesssim \tilde t_\xi |\xi|^{-\alpha} 
\lesssim |\xi|^{-1-|\alpha|} \big(\log(e+\tilde t_\xi)\big)^{2\nu}
\end{equation}
proves the desired result. \end{pf}

\begin{rem} If we diagonalise perfectly modulo $\H_{N,\nu}\{1\}$ we can construct the corresponding
fundamental solution $\mathcal E_\infty(t,s,\xi)=\widetilde{\mathcal E}_\infty(t,s,\xi)\mathcal Q_\infty(t,s,\xi)$ and the previous lemma holds for all multi-indices $\alpha$. Note, that $\H_{N,\nu}\{1\}$ is invariant under multiplications by $\widetilde{\mathcal E}_\infty(t,s,\xi)$.
\end{rem}

\subsubsection{Estimate in the oscillating subzone} It remains to estimate the fundamental solution 
in the oscillating sub-zone. In this part of the phase space it suffices to apply one step of diagonalisation,
i.e. diagonalisation modulo $\Sc_{N,\nu}^{\ell_1,\ell_2-2}\{-1,2\}$. We construct the fundamental solution $\mathcal E_1(t,s,\xi)=\widetilde{\mathcal E_1}(t,s,\xi)\mathcal Q_1(t,s,\xi)$ of $\D_t-\mathcal D-F_1-R_1$ following the lines of the previous subsection replacing the uniform integrability of the remainder by the estimate
\begin{multline}\label{eq:2.39}
  \int_{t_\xi}^{\tilde t_\xi} ||\mathcal R_1(t,t_\xi,\xi)||\d t 
  \lesssim \int_{t_\xi}^{\tilde t_\xi} 
  \frac{\big(\log(e+t)\big)^{2\nu}}{(1+t)^2|\xi|}
  \d t \\
  \sim 
  \frac{\big(\log(e+t)\big)^{2\nu}}
  {(1+t)|\xi|}\bigg|_{t_\xi}^{\tilde t_\xi} 
\le N \big(\log(e+t_\xi)\big)^\nu.
\end{multline}
This implies 
\begin{equation}\label{eq:2.36}
  ||\mathcal Q_1(t,s,\xi)|| \lesssim \exp\left( C'\big(\log(e+t_\xi)\big)^\nu   \right)
\end{equation}
for $t_\xi\le s\le t\le \tilde t_\xi$. For later use we need estimates for derivatives of $\mathcal Q_1$ with respect to $\xi$. 

\begin{lem}\label{lem:Q_Zosc}
  Assume (A1)$_{\ell_1,\ell_2}$, (A2) and (A3), with 
  $\ell_2\ge 2$. Then the matrix $\mathcal Q_1(t,s,\xi)$ satisfies
   \begin{equation}
     \|\D_\xi^\alpha \mathcal Q_1(t,t_\xi,\xi)\| \lesssim |\xi|^{-|\alpha|} (\log(e+ t))^{(1+2\nu)|\alpha|}  \exp\left( C'\big(\log(e+t_\xi)\big)^\nu   \right)
   \end{equation}
   for all $|\xi|\le N$ and all multi-indices $|\alpha|\le \ell_1 $.
\end{lem}

\begin{pf}
See, e.g., \cite[Lemma 2.9]{RS05a}. The proof follows directly by differentiating \eqref{eq:Qkrep} (for $k=1$), using the symbol estimates together with \eqref{eq:2.10} and estimates \eqref{eq:2.39}.
\end{pf}

\subsection{Combination of results}
We combine the representations of the two previous sections. We distinguish between three cases.
If $|\xi|>N$ then only the hyperbolic zone occurs. Hence
\begin{subequations}
\begin{equation}
\mathcal E(t,0,\xi)=M^{-1}(t,\xi)N_k^{-1}(t,\xi)\mathcal E_k(t,0,\xi)N_k(0,\xi)M(0,\xi).
\end{equation}
For $|\xi|\le N$ and $t\le t_\xi$ only the pseudo-differential zone is of interest and the fundamental 
solution is constructed in Section~\ref{sec2.1}. For $t\ge t_\xi$ we obtain
\begin{equation}
\mathcal E(t,0,\xi)=M^{-1}(t,\xi)N_k^{-1}(t,\xi)\mathcal E_k(t,t_\xi,\xi)N_k(t_\xi,\xi)M(t_\xi,\xi)
\mathcal E(t_\xi,0,\xi)
\end{equation}
\end{subequations}
with $k=1$ in the oscillating sub-zone and large $k$ and $\tilde t_\xi$ in place of $t_\xi$ in the regular sub-zone.

We will bring these representations into a unified form. For this,
let us define 
\begin{equation}\label{EQ:vala}
  \varphi_j(t,\xi)=\frac1t\int_0^t \lambda_j(\tau,\xi)\ \d\tau.
\end{equation}
Then the following statement holds true.
\begin{thm}\label{thm:repSol} 
Assume (A1)$_{\ell_1,\ell_2}$, (A2), (A3) and (A4). 
Then the fundamental matrix $\mathcal E(t,0,\xi)$ can be represented as
\begin{equation}
\mathcal E(t,0,\xi)=\sum_{j=1}^m e^{\i t |\xi| 
\varphi_j(t,\xi/|\xi|)} B_j(t,\xi),
\end{equation}
with matrices $B_j(t,\xi)\in\C^{m\times m}$ subject to the estimates
\begin{equation}
  \| B_j(t,\xi)\| \le  C  \exp\left(C'\big(\log(e+t)\big)^\nu\right)
\end{equation}
in $Z_{pd}(N,\nu)\cup Z_{osc}(N,\nu)$, and
\begin{equation}\label{EQ:reg-B-ests}
   \left\|\D_\xi^\alpha B_j(t,\xi)\right\| \le C_\alpha |\xi|^{-|\alpha|} (\log(e+ t))^{(2\nu+1)|\alpha|}
   \exp\left(C'\big(\log(e+t)\big)^\nu\right)
\end{equation}
for all $|\alpha|\le \min(\ell_1,\ell_2/2-1)$ in $Z_{reg}(N,\nu)$. 
\end{thm}
\begin{pf} 
Note first, that for $(t,\xi)\in Z_{pd}(N,\nu)\cup Z_{osc}(N,\nu)$,
all the terms $e^{\pm\i t|\xi|\phi_j(t,\xi/|\xi|)}$ are uniformly bounded and satisfy the corresponding symbolic estimates.
Thus artificially introducing these terms does not destroy our statement.

The estimate in the pseudo-differential zone and the oscillating sub-zone follows directly from Lemma~\ref{lemPD1} and equations \eqref{eq:2.26} and \eqref{eq:2.36}. Note that for $\nu=0$ we obtain just uniform bounds.

In the regular sub-zone the matrices $B_j$ collect terms from the diagonalisers $M$, $N_k$, the matrices $\mathcal E(t_\xi,0,\xi)$ and $\mathcal Q_k(t,\xi)$ and the terms arising from
\begin{equation}
\exp\left(\i\int_{\max(t_\xi,0)}^t F_k(\tau,\xi)\d\tau\right),
\end{equation}
all of which satisfy symbolic estimates due to 
Lemma~\ref{lemPD2}, \ref{lem:Zreg} and 
for the last one due to 
$F_k\in \Sc_{N,\nu}^{\ell_1,\ell_2-k}\{0,1\}$. Note further, that \eqref{eq:est-txi-tilde}
implies from $N_k(t,\xi)\in\Sc_{N,\nu}^{\ell_1,\ell_2-k}\{0,0\}$ 
\begin{equation}
  \| \D_\xi^\alpha N_k(\tilde t_\xi,\xi) \| \le C |\xi|^{-|\alpha|} \big(\log(e+\tilde t_\xi)\big)^{2\nu|\alpha|}
\end{equation}
which gives the $2\nu$ $\log$-estimate in the regular zone despite the better estimate of Lemma~\ref{eq:est-txi-tilde}.
\end{pf}

\section{Examples}\label{SEC:examples}
We will give some examples to illustrate the applicability of our previous construction. 

\subsection{Differential hyperbolic systems}
If we consider differential systems, they can always be written as $\D_t U = A(t,\D)U$, with 
\begin{equation}\label{EQ:dif-system}
    A(t,\D) = \sum_{j=1}^n A_j(t)\D_{x_j} + B(t),
\end{equation}
for $t$-dependent matrices $A_j(t), B(t)\in\C^{m\times m}$. Assumption (A1)$_{\ell_1,\ell_2}$ is satisfied if $A_j(t)\in \mathcal T_\nu\{0\}$ and
$B(t)\in\mathcal T_\nu\{1\}$, where 
\begin{equation}
   \mathcal T_\nu\{\rho\} = \left\{ f \in C^\infty(\R_+) : |\partial_t^k f(t)| \le C_k \left( \frac1{1+t} \big(\log(e+t)\big)^\nu \right)^{\rho-k}\right \},
\end{equation}
together with the assumption that the associated homogeneous matrix
\begin{equation}\label{eq:3.3}
   \sum_{j=1}^n A_j(t)\xi_j
\end{equation}
has real eigenvalues for all $t$ and $\xi$. The system is symmetric hyperbolic if the latter matrix is self-adjoint, which means that all matrices $A_j(t)$ are self-adjoint. We denote the eigenvalues  of \eqref{eq:3.3} by $\lambda_1(t,\xi),\dots,\lambda_m(t,\xi)$ in ascending order. Assumption (A2) guarantees that they are uniformly distinct if restricted to $\R\times\S^n$.

For the following we assume the system to be symmetric hyperbolic. Then Assumption (A3) simplifies to 
\begin{equation}
 \sup_{s,t} \left\|  \int_s^t \diag \left(M^{-1}(\theta,\xi)(\Im B(\theta))M(\theta,\xi) \right) \d\theta \right\| < \infty
\end{equation}
for the (0-homogeneous) unitary diagonaliser $M(t,\xi)$ of \eqref{eq:3.3}, while (A4) reads as
\begin{equation}
    \Im B(t) + \gamma(t) I \ge 0 
\end{equation}
for a suitable function $\gamma(t)$ satisfying \eqref{eq:hCond}. Both assumptions are trivially satisfied if we assume
$\Im B(t)\in L^1(\R_+)$. 

This corresponds to the situation in \cite{dAbbico:2009}. Main difference between the results of \cite{dAbbico:2009}
and our results is, that we do not subsequently require the non-vanishing of the Gaussian curvature of the 
characteristics. In Theorems \ref{THM:main} and \ref{THM:main-nonconvex} we derive 
estimates for solutions without this extra assumption. The case of \cite{dAbbico:2009} is covered by
the case $\gamma=2$ in Theorem \ref{THM:main}.

\subsection{Hyperbolic equations of second order} Second order equations under related assumptions are considered in \cite{RS05a}, \cite{Reissig2010} and related papers. If we consider the equation
\begin{equation}
  \D_t^2 u = \D_t\sum_{j=1}^n \D_{x_j} b_j(t) u + \sum_{1\le i\le j\le n} a_{i,j}(t) \D_{x_i}\D_{x_j} u + c(t) \D_t u + \sum_{j=1}^n d_j(t) \D_{x_j} u + e(t) u 
\end{equation}
with $a_{i,j}(t), b_j(t)\in\mathcal T_\nu\{0\}$, $c(t), d_j(t)\in\mathcal T_\nu\{1\}$ and $e(t)\in \mathcal T_\nu\{2\}$, we can rewrite it as $2\times 2$ pseudo-differential system in $U=(h(t,\D) u, \D_t u)^T$, $h(t,\xi) \simeq |\xi|_t$ (being $\xi$-independent within the pseudo-differential zone and smoothly continued to the hyperbolic zone by adding a suitable term from $\mathcal H_{N,\nu}\{1\}$),
\begin{equation}
   \D_t U = A(t,\D) U,   
\end{equation}
where within $Z_{hyp}(N,\nu)$ we have
\begin{multline}
  A(t,\xi) = \begin{pmatrix}
    0 & |\xi|\\
    |\xi|^{-1} \sum a_{i,j}(t) \xi_i\xi_j &   \sum b_j(t) \xi_j 
  \end{pmatrix}\\
  +
  \begin{pmatrix}
    0 & 0\\  |\xi|^{-1} \sum d_j(t) \xi_j & c(t)
  \end{pmatrix}
  +  \begin{pmatrix}
    0 & 0\\  |\xi|^{-1} e(t) & 0
  \end{pmatrix}
\end{multline}
modulo lower order terms form $\mathcal H_{N,\nu}\{1\}$. The second line gives the terms from $\Sc_{N,\nu}\{0,1\}$
and $\Sc_{N,\nu}\{-1,2\}$, respectively. Similarly, we obtain within $Z_{pd}(N,\nu)$ 
\begin{multline}
  A(t,\xi) = \begin{pmatrix}
  0& h(t)  \\
    h(t)^{-1} \sum a_{i,j}(t) \xi_i\xi_j &   \sum b_j(t) \xi_j 
  \end{pmatrix} \\
   +
  \begin{pmatrix}
      D_t \log h(t)  & 0\\  h(t)^{-1} \sum d_j(t) \xi_j & c(t)
  \end{pmatrix}
  +  \begin{pmatrix}
    0 & 0\\  h(t)^{-1} e(t) & 0
  \end{pmatrix}
\end{multline}
with $h(t)=h(t,0)=\frac N{1+t}(\log (e+t))^\nu$. Note that some of the `lower order entries' are the dominant terms now.

Assumptions (A1)$_{\ell_1,\ell_2}$ and (A2) are satisfied if we require that eigenvalues of the homogeneous principal part
\begin{equation}
  \lambda_\pm (t,\xi) = \frac12 \sum_{j=1}^n b_j(t)\xi_j \pm \sqrt{ \frac14 \bigg(\sum_{j=1}^n b_j(t)\xi_j\bigg)^2 +\sum_{i\le j} a_{i,j}(t) \xi_i \xi_j}
\end{equation}
are real and distinct, i.e., if $b_j(t)$ is real and 
\begin{equation}
0 < \bigg(\sum_{j=1}^n b_j(t)\xi_j\bigg)^2 + 4 \sum_{i\le j} a_{i,j}(t) \xi_i \xi_j, \qquad \xi\ne0.
\end{equation}

Assumption (A3) requires that the integral of
\begin{multline}
   \frac{\lambda_{\pm} (t,\xi)}{\lambda_+(t,\xi)-\lambda_-(t,\xi)} \Im  c(t) +
   \frac{\partial_t\lambda_\pm(t,\xi)}{\lambda_+(t,\xi)-\lambda_-(t,\xi)} +  \frac{ \Im \sum_{i} d_i(t)\xi_i}{\lambda_+(t,\xi)-\lambda_-(t,\xi)}
\end{multline}
is uniformly bounded over the hyperbolic zone. This is in particular the case if $\Im c(t)$ and $\Im d_i(t)$ are integrable and the equivalent to assumption \eqref{EQ:hh-int} is satisfied. 
Furthermore,  (A4) becomes a condition on the lower order terms guaranteeing that their influence is dominated by the principle part. 

\subsection{Hyperbolic equations of higher order}\label{expl:3} As third example 
we consider hyperbolic equations of higher order,
\begin{equation}\label{eq:hh-CP}
   \D_t^m u + \sum_{j=0}^{m-1} 
   \sum_{j+|\alpha|\le m} a_{j,\alpha}(t) \D_t^j \D_x^\alpha u = 0,
\end{equation}
together with the Cauchy data
\begin{equation}
   \D_t^j u(0,\cdot) = u_j, \qquad j=0,1,\ldots, m-1.
\end{equation}
Similar to the second order case, this can be brought into pseudo-differential system form.

Then (A1)$_{\ell_1,\ell_2}$ is satisfied if we assume $a_{j,\alpha}\in \mathcal T_\nu\{m-j-|\alpha|\}$ and if the roots $\tau=\lambda_j(t,\xi)$ of the homogeneous principal part, i.e., the solutions of 
\begin{equation}
   \tau^m +\sum_{j=0}^{m-1} \sum_{j+|\alpha|=m} a_{j,\alpha}(t) \tau^j \xi^\alpha = 0,
\end{equation}
are real. For (A2) we require that they are distinct. If we consider equations which are homogeneous of order $m$, i.e., $ a_{j,\alpha}(t)=0$ for $j+|\alpha|<m$, condition (A3)  is equivalent to
\begin{equation}\label{EQ:hh-int}
   \max_{1\le j\le m} \sup_{T\ge 0}\sup_{\xi\ne0} \left| \sum_{k\ne j} \int_0^T \frac{\partial_t \lambda_j(t,\xi)}{\lambda_j(t,\xi)-\lambda_k(t,\xi)}\d t  \right|<\infty.
\end{equation}
while  Assumption (A4) is trivially satisfied. 

\begin{rem}\label{REM:MR}
In \cite{MR10} homogeneous equations of order $m$ where considered 
under the stronger assumption 
$\partial_t a_{j,\alpha}(t)\in L^1(\R_+)$ for $j+|\alpha|=m$.
This implies, in particular, that \eqref{EQ:hh-int} is 
satisfied. The representation of solutions is obtained there
by the asymptotic integration method, requiring less conditions
on the coefficients $a_{j,\alpha}$ (conditions on only the first 
time derivatives are enough). However, the asymptotic integration method
yields a representation of solutions with less control
on its amplitudes and an additional loss of regularity in the
dispersive estimates, which does not occur with the
present method.
\end{rem}

\section{Dispersive estimates}

In this section we are concerned with dispersive estimate for solutions represented in the form of 
Theorem~\ref{thm:repSol}. From now on we assume that $n\ge2$.
 
First, we give the estimate for low frequencies:

\begin{lem}\label{LEM:low-freq}
Assume (A4). Then solution $U=U(t,x)$ 
to \eqref{eq:CP} satisfies
\begin{multline}
   \| \Fcal^{-1}[(1-\chi_{reg}(t,\xi))\widehat  U(t,\xi) ] 
   \|_{L^\infty(\R^n)} \\   \le C (1+t)^{-n} \big(\log(e+t)\big)^{2n\nu} \exp(C'\log(e+t)^{\nu})   
   \|U_0\|_{L^1(\R^n)}
\end{multline}
localised to $Z_{pd}(N,\nu)\cup Z_{osc}(N,\nu)$ 
(for any choice of $N$).
\end{lem}
\begin{pf}
Based on the mapping property of the Fourier transform $\Fcal : L^1(\R^n)\to L^\infty(\R^n)$ 
and H\"older inequality it is sufficient to estimate  
$$\| \mathcal E(t,0,\xi)\chi_{pd+osc}(t,\xi)\|_{L^1(\R^n)}\le\| 
\mathcal E(t,0,\xi)\|_{L^\infty(|\xi|\le \tilde \xi_t)} 
\|\chi_{pd+osc}\|_{L^1(\R^n)}$$ and, therefore, the estimate follows 
from Lemma~\ref{lemPD1} and the geometry of the zone.
\end{pf}

For the high-frequency part we at first derive some abstract statements giving $L^p$--$L^q$ decay rates for 
oscillatory integrals with a related structure to the ones constructed in Theorem~\ref{thm:repSol},
$$
   T_t f(x) = \int_{\R^n} 
   \mathrm e^{\i(x\cdot\xi + t \varphi(t,\xi)} b(t,\xi) \widehat f(\xi) \d\xi,
$$
with a real phase function $\varphi$ and and amplitude $b$. 
For simplicity, we omit the inclusion of logarithmic terms in 
assumptions here, later they will give a (for $\nu<1$ small) change in the 
decay rates by a simple argument. 

We introduce a cut-off function of the form
$\psi((1+|t|)\xi)$ for some $\psi\in C_0^\infty(\Rn)$
such that $\psi(\xi)\equiv 1$ for 
$\vert \xi \vert \le \frac12$, 
and $0$ for  $\vert \xi \vert \ge 1$. 
We recall that we use the notation 
$\dot{L}^{p}_{\kappa}(\Rn)$ for the
homogeneous Sobolev space $\dot{W}^{\kappa}_p(\Rn).$

For now we may assume that phase functions
 satisfy $\varphi(t,\xi)>0$ uniformly for all 
$\xi\ne0$ in the support of $b(t,\xi)$. This can always be achieved by localisation and in combination with adding a linear function to the phase which will not affect estimates. 
Dispersive estimates for the corresponding Fourier 
integral operator $T_t$ are determined by the geometric
properties of the family of  Fresnel surfaces  
\begin{equation}
  \Sigma^t = \{\, \xi\in\R^n\setminus\{0\}\,:\, 
  \varphi(t,\xi) = |\xi| \varphi(t,\xi/|\xi|) =1\, \}.
\end{equation}
Following \cite{Sug96} we introduce two indices 
for such a Fresnel surface $\Sigma$, assuming that
it is of class $C^k$ with $k$ being sufficiently large. 
For $u\in\Sigma$ we denote by $\mathcal T_u$ the  
tangent hyperplane to $\Sigma$ at $u$. Then for any plane $H$ 
of dimension 2 which contains $u$ and the normal of
$\Sigma$ at $u$ we denote by $\gamma(\Sigma;u,H)$ 
the order of contact between the curve
$\Sigma\cap H$ and its tangent $H\cap \mathcal T_u$. Furthermore, we set
\begin{equation}
 \gamma(\Sigma) = \sup_{u} \sup_{H} \gamma(\Sigma;u,H),\qquad 
 \gamma_0 (\Sigma) = \sup_u \inf_H \gamma(\Sigma; u,H).
\end{equation}
Obviously the definition implies 
$2\le\gamma_0(\Sigma)\le\gamma(\Sigma)$. 
For isotropic problems $\Sigma$ is a dilation of $\S^{n-1}$ and 
$\gamma(\S^{n-1})=\gamma_0(\S^{n-1})=2$.
Moreover, if the Gaussian curvature of $\Sigma$ never
vanishes we have $\gamma(\Sigma)=2$.

To control the order of contact of $\Sigma$ by tangent lines
quantitatively, we will introduce another quantity
$\varkappa(\Sigma)$. First, assume that $\Sigma$ is convex and
that it is of class $C^{\gamma(\Sigma)+1}$.
For $u\in\Sigma$, rotating $\Sigma$ if
necessary, we may assume that it is parameterised
by points $\{(y,h(y)), y\in \Omega\}$ near $u$ for an open set
$\Omega\subset\R^{n-1}$. For $u=(y,h(y))$, let us define
$$\varkappa(\Sigma;u)=\inf_{|\omega|=1} \sum_{j=2}^{\gamma(\Sigma)}
\left| \frac{\partial^j}{\partial\rho^j} h(y+\rho\omega)|_{\rho=0}
\right|.$$ 
From the definition of $\gamma(\Sigma)$
it follows that $\varkappa(\Sigma;u)>0$ for all $u\in\Omega$.
Indeed, from the definition of $\gamma(\Sigma;u,H)$ it follows
that if $\omega$ is such that $y+\rho\omega\in H$, then
$$
\varkappa(\Sigma;u,H)=\left|\frac{\partial^{\gamma(\Sigma;u,H)}}
{\partial\rho^{\gamma(\Sigma;u,H)}} h(y+\rho\omega)|_{\rho=0}
\right|>0.
$$
Now, we clearly have
$\sum_{j=2}^{\gamma(\Sigma)}
\left| \frac{\partial^j}{\partial\rho^j} h(y+\rho\omega)|_{\rho=0}
\right|\geq \varkappa(\Sigma;u,H)$,
and hence we have $\varkappa(\Sigma;u)>0$ since the set
$|\omega|=1$ is compact.
Noticing that $\varkappa(\Sigma;u)$
is a continuous function of $u$, by compactness of
$\Sigma$ it
follows that if we define
$$
\varkappa(\Sigma)=\min_{u\in\Sigma} \varkappa(\Sigma;u),
$$
then $\varkappa(\Sigma)>0$.

If $\Sigma$ is not convex, we define
$$
\varkappa_0(\Sigma)=\min_{u\in\Sigma} 
\sup_{|\omega|=1} \sum_{j=2}^{\gamma_0(\Sigma)}
\left| \frac{\partial^j}{\partial\rho^j} h(y+\rho\omega)|_{\rho=0}
\right|.
$$
Again, we have $\varkappa_0(\Sigma)>0$.

The quantities $\varkappa(\Sigma^t)$ and $\varkappa_0(\Sigma^t)$
evaluated for all times $t$ will allow us to ensure
that if the surfaces degenerate, it will be at least
uniformly with respect to $t$. If the degeneracy of a certain
order is not uniform, it may become uniform when we increase
integers $\gamma$ in the formulations below.

For completeness we include the following proposition.  It gives a 
criterion on the convexity of the level sets of the phases,
as well as an upper bound on the index $\gamma$ for systems arising from differential operators.
\begin{prop}\label{PROP:positivity} Assume that $\det A_1(t,\xi)$ is polynomial in $\xi$ and denote 
by  $\lambda_k(t,\xi)$, $k=1,\ldots,m$, the characteristic roots of the
operator \eqref{eq:CP}, ordered by
$\lambda_1(t,\xi)<\lambda_2(t,\xi)<\cdots<\lambda_m(t,\xi)$ for $\xi\not=0$.
Suppose that all the Hessians $\nabla_\xi^2
\lambda_k(t,\xi)$ are
semi-definite for $\xi\not=0$. Then there exists a polynomial
$\alpha(t,\xi)$ in $\xi$, of order one, such that
$\lambda_{m/2}(t,\xi)<\alpha(t,\xi)<\lambda_{m/2+1}(t,\xi)$ 
{\rm (}if $m$ is even{\rm )} or
$\alpha(t,\xi)=\lambda_{(m+1)/2}(t,\xi)$ 
{\rm (}if $m$ is odd{\rm )}. 
Moreover\footnote{Here, as in \eqref{EQ:vala},
we define
$\varphi_j(t,\xi)=\frac1t\int_0^t \lambda_j(\tau,\xi)\ \d\tau.$},
the hypersurfaces $\Sigma_k^t=\{\xi\in\Rn;\, 
\widetilde{\va}_k(t,\xi)=\pm 1\}$
with $\widetilde{\va}_k(t,\xi)=\va_k(t,\xi)-\alpha(t,\xi)$ 
$(k\not=(m+1)/2)$ are convex and we have
$\gamma(\Sigma_k^t)\leq 2\lfloor m/2\rfloor.$
\label{PROP:limitingphases}
\end{prop}
 By taking the determinant
 of the system \eqref{eq:CP} we can reduce the analysis of
 characteristics to scalar equations. Then 
this proposition can be readily shown by a modification 
of the argument in \cite{sugi94}, so we omit the details.

We are now ready to estimate the appearing oscillatory
integrals. In the sequel, for $r>0$, by $\lfloor r\rfloor$ we denote
its integer part. The Proposition 
below extends Proposition 4.2 in
\cite{MR10} where the limit of $\varphi(t,\xi)$ was
assumed to exist as $t\to\infty$, and which deals
with amplitudes of type $(0,0)$.

\begin{prop} \label{PROP:oscillatory-convex}
Let $\gamma\in\N$.
Let $T_t$ be an operator defined by
\[
T_t f(x)=\int_\Rn \mathrm{e}^{\i(x\cdot\xi+t\varphi(t,\xi))}
\ \br{1-\psi((1+|t|)\xi)}\
b(t,\xi)\ \widehat{f}(\xi)\ \d\xi,
\]
where $\varphi(t,\xi)$ is real valued, continuous in $t$, 
smooth in
$\xi\in\Rn\backslash 0$, homogeneous of order one in $\xi$, 
and such that for some $t_0>0$ and $C>0$ we have
\begin{equation}\label{EQ:ass-phi}
C^{-1}|\xi| \leq \varphi(t,\xi) \leq C|\xi| 
\;\;\text{and}
\;\; |\partial_\xi^\alpha\varphi(t,\xi)|\leq C|\xi|^{1-|\alpha|}
\;\;\text{for all}
\;\; t\geq t_0,\; \xi\not=0,
\end{equation} 
and all $\alpha$ such that $|\alpha|\leq\max\{\gamma+1,
\lfloor(n-1)/\gamma\rfloor+2\}.$
Assume that the sets
\begin{equation}\label{EQ:ass-sigma}
\Sigma^t=\{\xi\in\Rn\backslash 0:\varphi(t,\xi)=1\}
\end{equation}  
are convex for all $t\geq t_0$, and assume that 
$\lim\sup_{t\to \infty} \gamma(\Sigma^t)\leq\gamma$ and
that $\lim\inf_{t\to\infty} \varkappa(\Sigma^t)>0$.  
Let us suppose that the amplitude
$b(t,\xi)$ satisfies 
\begin{equation}\label{EQ:ass-amp}
|\partial_\xi^\alpha b(t,\xi)|\leq C_\alpha |\xi|^{-|\alpha|}
\quad\text{for all}
\quad |\alpha|\leq \lfloor(n-1)/\gamma\rfloor+1.
\end{equation} 
Let $1< p\leq 2\leq q< \infty$ be such that
$\frac{1}{p}+\frac{1}{q}=1$.
Then for all $t\geq t_0$ we have the estimate
\[
 \Vert T_t f \Vert_{L^q(\Rn)} \leq 
 C t^{-\frac{n-1}{\gamma}\p{\frac{1}{p}-
 \frac{1}{q}}}
 \Vert f \Vert_{\dot{L}^p_{N_p}(\Rn)},
\]
where $N_p=\left(n-\frac{n-1}{\gamma}
\right)
\left(\frac{1}{p}-\frac{1}{q}\right)$. 
\end{prop}

\begin{pf}
Since the $L^2$--estimate
$\Vert T_t f \Vert_{L^2}\leq C \Vert f \Vert_{L^2}$
readily follows from the Plancherel identity, by interpolation
we only need to prove that
\begin{equation}\label{EQ:est2}
\Vert \widetilde{T_t} f \Vert_{L^\infty(\Rn)} 
\leq C t^{-\frac{n-1}{\gamma}}
 \Vert f \Vert_{L^1(\Rn)},
\end{equation}
with $\widetilde{T_t}=T_t\circ |D|^{-N_1}$,
where the amplitude $a(t,\xi)$ of $\widetilde{T_t}$ satisfies
$|\partial_\xi^\alpha a(t,\xi)|\leq C_\alpha 
|\xi|^{-N_1-|\alpha|}$ for
all $|\alpha|\leq \lfloor(n-1)/\gamma\rfloor+1$ and
$N_1=n-\frac{n-1}{\gamma}$. Let us make an additional
decomposition
\begin{equation}\label{EQ:T-decomp}
\widetilde{T_t}=\widetilde{T_t^{(1)}}+\widetilde{T_t^{(2)}}
=\widetilde{T_t}\circ (1-\psi(D))+
\widetilde{T_t}\circ \psi(D).
\end{equation}
First we will treat the operator $\widetilde{T_t^{(1)}}$.
To simplify the notation, for this part we will denote
$a(t,\xi)(1-\psi(\xi))$ by $a(t,\xi)$ again.
By using Besov spaces, we can microlocalise the desired estimate
\eqref{EQ:est2}
to spherical layers in the frequency space. Indeed, let 
$\{\Phi_j\}_{j=0}^\infty$ be the Littlewood-Paley partition of
unity, and let 
$$
\Vert u \Vert_{B^s_{p,q}}=\p{ \sum_{j=0}^\infty \p{2^{js} \Vert \Fcal^{-1}
\Phi_j(\xi)\Fcal u \Vert_{L^p(\Rn)}}^q}^{1/q}
$$
be the norm of the Besov space $B^s_{p,q}$. Then, because of the
continuous embeddings $L^p\subset B^0_{p,2}$ for $1<p\leq 2$,
and $B^0_{q,2}\subset L^{q}$ for $2\leq q<+\infty$ 
(see \cite{BL}), it is sufficient to prove the uniform estimate
for the operators with amplitudes $a(t,\xi)\Phi_j(\xi)$.
Writing 
$\Phi_j(\xi)=\Phi_j(\xi)\Psi\p{\frac{
\varphi(t,\xi)}{2^j}}$ with
some function $\Psi\in C_0^\infty(0,\infty)$, we may prove
the uniform estimate for operators with amplitudes
$a(t,\xi)\Psi\p{\frac{
\varphi(t,\xi)}{2^j}}.$
Such choice of $\Psi$ is possible due to our assumption
\eqref{EQ:ass-phi} on $\varphi(t,\xi)$,
and we restrict the analysis for large enough $t$.
Let 
\[ 
I(t,x)=\int_\Rn \mathrm{e}^{i(x\cdot\xi+t\varphi(t,\xi))}
a(t,\xi)\Psi\p{\frac{\varphi(t,\xi)}{2^j}}\, \d\xi
\]
be the kernel of the corresponding operator.
Since we easily have the $L^2$--$L^2$ estimate by the 
Plancherel identity, by analytic interpolation
we only need to prove the $L^1$--$L^\infty$ case of \eqref{EQ:est2}.
In turn, this follows from the
estimate $|I(t,x)|\leq C t^{-\frac{n-1}{\gamma}},$
with constant $C$ independent of $j$.

Let $\kappa\in C_0^\infty(\Rn)$ be supported in 
a ball  centred at the origin, with some radius $r>0$ to be chosen
later. 
We decompose the kernel $I(t,x)$ as
\begin{align}\label{EQ:I-dec}
I(t,x) & =  I_1 (t,x)+I_2 (t,x) \\ \nonumber & =
\int_\Rn \mathrm{e}^{i(x\cdot\xi+t\varphi(t,\xi))}
a(t,\xi) \kappa\left(t^{-1}x+\nabla_\xi\varphi(t,\xi)\right)
\Psi\p{\frac{\varphi(t,\xi)}{2^j}} \d\xi
\\ \nonumber
& + 
\int_\Rn \mathrm{e}^{i(x\cdot\xi+t\varphi(t,\xi))}
a(t,\xi) (1-\kappa)
\left(t^{-1}x+\nabla_\xi\varphi(t,\xi)\right)
\Psi\p{\frac{\varphi(t,\xi)}{2^j}} \d\xi.
\end{align}
We can easily see that 
$|I_2(t,x)|\leq C t^{-\frac{n-1}{\gamma}}$. In fact,
we can show
$|I_2(t,x)|\leq C t^{-l}$ 
for $l=\lfloor(n-1)/\gamma\rfloor+1$ and then the required estimate
simply follows since $l>(n-1)/\gamma$. 
Indeed, on the support of $1-\kappa$,
we have $|x+t\nabla_\xi\varphi (t,\xi)|\geq rt>0$. Thus,
integrating by parts with operator
$P=\frac{x+t\nabla_\xi \varphi (t,\xi)}
{i|x+t\nabla_\xi\varphi (t,\xi)|^2}
\cdot \nabla_\xi$, we get
\begin{multline}\label{EQ:I2tx}
I_2(t,x)= \\
\int_\Rn 
\mathrm{e}^{i(x\cdot\xi+t\varphi(t,\xi))}
(P^*)^l \left[ a(t,\xi) (1-\kappa)
\left(t^{-1}x+\nabla_\xi\varphi(t,\xi)\right)
\Psi\p{\frac{\varphi(t,\xi)}{2^j}}\right]  \d\xi.
\end{multline}
Using the fact that $|\partial_\xi^\alpha\varphi(t,\xi)|
\leq C|\xi|^{1-|\alpha|}$ by \eqref{EQ:ass-phi}, 
we readily observe from \eqref{EQ:I2tx} that 
the required estimate $|I_2(t,x)|\leq C t^{-l}$ holds. Here we
also used the condition \eqref{EQ:ass-amp} assuring that
we can perform the integration by parts $\lfloor(n-1)/\gamma\rfloor+1$
times. We note that since there is one more $\nabla_\xi\va$
involved here, the condition $|\alpha|\leq \lfloor(n-1)/\gamma\rfloor+2$
in \eqref{EQ:ass-phi} allows us to integrate by parts $l$
times.

Now we will turn to estimating $I_1(t,x)$. 
Here we are going to use the structure of the sets
$\Sigma^t$ in \eqref{EQ:ass-sigma}, restricting to $t$ large
enough.
We recall that \eqref{EQ:ass-phi} implies, in particular, that 
$\va(t,\xi)>0$ for all $\xi\not=0$.
By rotation, we can always microlocalise in some narrow cone
around $e_n=(0,\ldots,0,1)$, and in this cone we can 
parameterise 
$$\Sigma^t=\{(y,h_t(y)): y\in U\}$$ 
for some open $U\subset\R^{n-1}$.
In other words, we have 
${\varphi}(t;y,h_t(y))=1$, and it follows that
$h_t$ is smooth and
$\nabla h_t:U\to \nabla h_t(U)\subset\R^{n-1}$ is a
homeomorphism. The function $h_t$ is concave if $\Sigma^t$ is convex.
We claim that 
\begin{equation}\label{EQ:ht}
|\partial_y^\alpha h_t(y)|\leq C_\alpha,\quad
\textrm{for all}\quad y\in U\quad \textrm{and large enough}
\quad t,
\end{equation}
for all $|\alpha|\leq\max\{\gamma+1,\lfloor(n-1)/\gamma\rfloor+2\}.$
Indeed, let us look at $|\alpha|=1$ first.
From ${\varphi}(t;y,h_t(y))=1$ we get
$\nabla_y {\varphi}+\partial_{\xi_n}
{\varphi}\cdot \nabla h_t(y)=0$.
From \eqref{EQ:ass-phi} we have 
$|\nabla_\xi{\varphi}|\leq C$,
so also $|\nabla_y {\varphi}|\leq C$.
By Euler's identity we have 
$\partial_{\xi_n} {\varphi} (t;e_n)=
{\varphi}(t;e_n)>0$,
so we have $|\partial_{\xi_n} {\varphi}|\geq c>0$
since we are in a narrow cone around $e_n$. From this it
follows that $|\nabla_y h_t(y)|\leq C$ for all $y\in U$
and $t$ large enough. A similar argument proves the boundedness
of higher order derivatives in \eqref{EQ:ht}.

Now, let us turn to analyse the structure of the sets $\Sigma^t$.
We have the Gauss map
$$
G:\Sigma^t\ni\zeta \mapsto \frac{
\nabla_\zeta{\varphi}(t;\zeta)}
{|\nabla_\zeta{\varphi}(t;\zeta)|}\in\Snm,
$$
and for $x=(x^\prime,x_n)\in\R^{n-1}\times\R$
near the point $-\nabla_\zeta{\varphi}(t;e_n)$
we define $z_t\in U$ by
$(z_t,h_t(z_t))=G^{-1}(-x/|x|)$. 
Then $(-\nabla_y h_t(y),1)$ is normal to $\Sigma^t$ at
$(y,h_t(y))$, so we get
$$
-\frac{x}{|x|}=\frac{(-\nabla_y h_t(z_t),1)}
{|(-\nabla_y h_t(z_t),1)|} \quad \textrm{and} \quad
\frac{x^\prime}{x_n}=-\nabla_y h_t(z_t).
$$
Making change of variables $\xi=(\widetilde\lambda y,
\widetilde\lambda h_t(y))$
and using ${\varphi}(t,\xi)=\widetilde\lambda$, we get
\begin{eqnarray*}
&& I_1(t,x)  \\
& =&
\int_0^\infty \int_U \mathrm{e}^
{i\widetilde\lambda(x^\prime\cdot y+x_n h_t(y)+t)}
a(t,\widetilde\lambda y,\widetilde\lambda h_t(y))
\Psi\p{\frac{\widetilde\lambda}{2^j}}\kappa_0(t,x,y) 
\abs{\frac{d\xi}{d(\widetilde\lambda,y)}} \, \d y\ \d\widetilde\lambda 
\nonumber \\
& = &
 \int_0^\infty \int_U \mathrm{e}^
{i\widetilde\lambda(-x_n\nabla_y h_t(z_t)\cdot y+x_n h_t(y)+t)}
\left[\widetilde\lambda^l a(t,\widetilde\lambda y,
\widetilde\lambda h_t(y))\right]
\widetilde\lambda^{n-1-l} \times \nonumber \\
& &\quad \times
\Psi\p{\frac{\widetilde\lambda}{2^j}}\kappa_0(t,x,y) \chi(t,y) 
\, \d y\ \d\widetilde\lambda \nonumber\\
& = & 
\int_0^\infty \int_U \mathrm{e}^
{i\lambda(-\nabla_y h_t(z_t)\cdot y+ h_t(y)+t x_n^{-1})}
\widetilde{a}(t,x_n,\lambda y,\lambda h_t(y)) 
\lambda^{n-1-l}  \times \nonumber \\
& & \quad \times \Psi\p{\frac{\lambda}{2^j x_n }} x_n^{-n+1+l-1}
\kappa_0(t,x,y) \chi(t,y) \, \d y\ \d\lambda,\nonumber 
\end{eqnarray*}
where $\kappa_0(t,x,y)=
\kappa\left(t^{-1}x+\nabla_\xi\va(t;y,h_t(y))\right)$,
and 
$$\widetilde{a}(t,x_n,\lambda y,\lambda h_t(y))=
(x_n^{-1}\lambda)^l a
\left(t,x_n^{-1} \lambda y,x_n^{-1} \lambda h_t(y)\right),$$
and
where we made a change $\widetilde\lambda=x_n^{-1}\lambda$ in the
last equality. Here also we used
$\big|\frac{\d\xi}{\d (\widetilde\lambda,y)}\big|=
\widetilde\lambda^{n-1} \chi(t,y)$,
where $\chi(t,y)$ and its derivatives with respect to $y$
are bounded because of \eqref{EQ:ht}.

If we choose $r$ in the definition of the cut-off function
$\kappa$ sufficiently small, then on its support we have 
$|x|\approx |x_n| \approx t$, and we can estimate
\begin{equation}\label{EQ:IJ}
\begin{aligned}
|I_1(t,x)| & \leq C t^{-n+l}\int_0^\infty 
\abs{J(\lambda,z_t)
\Psi\p{\frac{\lambda}{2^j t }}
\lambda^{n-1-l}} \, \d\lambda \\ 
& =
C t^{-n+l}2^{j(n-l)}\int_0^\infty 
\abs{J(2^j \lambda,z_t)
\Psi\p{\frac{\lambda}{t }}
\lambda^{n-1-l}} \, \d\lambda,
\end{aligned}
\end{equation}
with
$$
J(\lambda,z_t)=\int_U \mathrm{e}^
{i\lambda(-\nabla_y h_t(z_t)\cdot y+ h_t(y)+t x_n^{-1})}
\widetilde{a}(t,x_n,\lambda y,\lambda h_t(y)) 
\kappa_0(t,x,y) \chi(t,y)\ \d y.
$$
We will show that
\begin{equation}\label{EQ:estJ}
\abs{J(\lambda,z_t)}\leq C(1+\lambda)^{-\frac{n-1}{\gamma}},
\quad \lambda>0.
\end{equation}
Then, if we take $l=n-\frac{n-1}{\gamma}$, and use
\eqref{EQ:IJ} and \eqref{EQ:estJ}, we get
\begin{equation}\label{EQ:I1tx-est2}
\begin{aligned}
|I_1(t,x)| & \leq C t^{-\frac{n-1}{\gamma}} 
2^{j\frac{n-1}{\gamma}}
\int_0^\infty
(2^j\lambda)^{-\frac{n-1}{\gamma}}
\Psi\p{\frac{\lambda}{t}}\lambda^{\frac{n-1}{\gamma}-1} \, \d\lambda
\\
& = C t^{-\frac{n-1}{\gamma}} 
\int_0^\infty
\lambda^{-1}
\Psi\p{\frac{\lambda}{t}} \, \d\lambda 
=  C t^{-\frac{n-1}{\gamma}} 
\int_0^\infty
\lambda^{-1}
\Psi\p{\lambda}\, \d\lambda \\
& \leq C t^{-\frac{n-1}{\gamma}},
\end{aligned}
\end{equation}
which is the desired estimate
for $I_1(t,x)$.

Let us now prove \eqref{EQ:estJ}. We note first that
with this choice of $l$ we have
\begin{equation}\label{eq:a-tilde}
 |\partial^\alpha_y \widetilde{a}|\leq C \quad
 \textrm{for all}\quad |\alpha|\leq 
 \left\lfloor (n-1)/\gamma \right\rfloor+1.
\end{equation}
Now, estimate \eqref{EQ:estJ} follows from 
Theorem \ref{THM:oscintthm} in Appendix A. Indeed, we can write
$J(\lambda,z_t)$ in polar coordinates
$(\rho,\omega)$ with $y=\rho\omega+z_t$, so that
\begin{equation}\label{EQ:gather1}
J(\lambda,z_t)=
e^{i\lambda(t x_n^{-1}+h_t(z_t))}
\int_{{\mathbb S}^{n-2}} \int^{\infty}_{0}
e^{i\lambda F(\rho,z_t,\omega)} \beta(\rho,z_t,\omega) 
\rho^{n-2} \ \d\rho\ \d\omega,
\end{equation}
with
\begin{gather}
\label{osc-phase}
F(\rho,z_t,\omega)= h_t(\rho\omega+z_t)-h_t(z_t)-
\rho \nabla_y h_t(z_t)\cdot\omega, \\
\nonumber
\beta(\rho,z_t,\omega)=
\widetilde{a}\left(t,x_n,\lambda (\rho\omega+z_t),\lambda h_t
(\rho\omega+z_t)\right) 
\kappa_0(t,x,\rho\omega+z_t) \chi(t,\rho\omega+z_t),
\end{gather}
where we can assume in addition that $\chi=0$ unless
$\rho\omega+z_t\in U$, so both $\rho$ and $\omega$ vary over
bounded sets.
Now, we want to apply Theorem \ref{THM:oscintthm}
to obtain estimate \eqref{EQ:estJ}.
The function $F$ in \eqref{osc-phase} satisfies
condition (F3) of Theorem \ref{THM:oscintthm} 
because of the definition of the
convex index $\gamma$ and because $h_t$ is concave.
From the assumption 
$\lim\inf_{t\to\infty} \varkappa(\Sigma_t)>0$ it follows that 
function $F$ satisfies property (F2) of
Theorem \ref{THM:oscintthm}. One can readily see
that the other conditions of Theorem \ref{THM:oscintthm}
are satisfied, implying \eqref{EQ:estJ}.

We now prove the estimates for the integral
$\widetilde{T_t^{(2)}}$ in the decomposition \eqref{EQ:T-decomp}.
This proof is essentially similar to 
that for $\widetilde{T_t^{(1)}}$ with a few differences
that we will point out here. Again, by interpolation,
it is sufficient to prove estimate
\begin{equation}\label{EQ:est-Tt2}
 \Vert \widetilde{T_t^{(2)}} f \Vert_{L^\infty(\Rn)} \leq 
 C  t^{-\frac{n-1}{\gamma}}
 \Vert f \Vert_{{L}^1(\Rn)},
\end{equation} 
with amplitude $a(t,\xi)$ satisfying
\begin{equation} \label{EQ:assAMP-low-freq1}
|\partial_\xi^\alpha a(t,\xi)|\leq C_\alpha |\xi|^{-N_1-|\alpha|}
\quad\text{for all}
\quad |\alpha|\leq \lfloor(n-1)/\gamma\rfloor+1,
\end{equation}
with $N_1=n-\frac{n-1}{\gamma}$.

We continue to estimate the integrals 
$I_1(t,x)$ and $I_2(t,x)$ in the decomposition
\eqref{EQ:I-dec} corresponding to the integral
$\widetilde{T_t^{(2)}}$.
In general, since for $\widetilde{T_t^{(2)}}$
we work with low frequencies $|\xi|<1$ only,
no Besov space decomposition is necessary, so we do not need
to introduce function $\Psi$ and $\Phi_j$, and we can take
$\Psi=1$.

The additional complications are related to the fact that 
in principle derivatives of the amplitude of the operator
$\widetilde{T_t^{(2)}}$ may introduce
an additional growth with respect to $t$.
In the estimate for $I_2(t,x)$ we performed integration
by parts with operator $P$. Now after the integration by parts
the amplitude of this integral in \eqref{EQ:I2tx} is
$$
(P^*)^l\br{\br{1-\psi((1+|t|)\xi)} a(t,\xi) (1-\kappa)
\left(t^{-1}x+\nabla_\xi\varphi(t,\xi)\right)
\Psi\p{\frac{\varphi(t,\xi)}{2^j}}}.
$$
Then, if any of the $\xi$-derivatives falls on 
$\br{1-\psi((1+|t|)\xi)}$, we get an extra factor $t$ which
is cancelled with $t^{-1}$ in the definition of $P$.
However, in this case we can then restrict to the support of
$\nabla\psi$ which is contained in the ball with radius
$(1+|t|)^{-1}$, so we are in the situation of
very low frequencies $|\xi|\leq C t^{-1}$ again. 
Consequently, for this integral we actually get a 
better decay rate of Lemma \ref{LEM:low-freq}.
If none of the derivatives in $(P^*)^l$ fall on
$\br{1-\psi((1+|t|)\xi)}$, the argument is the same as
in the proof of the estimate for $I_2(t,x)$ for the
integral $\widetilde{T_t^{(1)}}$.

Let us now analyse the term 
$I_1(t,x)$ corresponding to
$\widetilde{T_t^{(2)}}$. Recall now that in the 
process of estimating $I_1(t,x)$ corresponding to
$\widetilde{T_t^{(1)}}$, 
we made a change of variables $\widetilde\lambda=x_n^{-1}\lambda$.
As it was then pointed out, if $r$ in the definition of
the cut-off function $\kappa$ is chosen sufficiently small,
on its support we have $|x_n|\approx |t|$. On the
other hand, we have $|\xi|\approx \widetilde\lambda$ 
by the definition
of $\widetilde\lambda$, since we may assume that 
$\varphi(t,\xi)$ is strictly positive for $\xi\not=0$. 
It then follows
that $(1+|t|)\xi\approx \widetilde\lambda |x_n|\approx\lambda$, 
and so the change of variables $\widetilde\lambda=x_n^{-1}\lambda$
changes $\br{1-\psi((1+|t|)\xi)}$ into 
$\br{1-\widetilde{\psi}(\lambda)}$ in the amplitude of $I_1(t,x)$.
Justifying this argument, we can then continue as in
the case of $\widetilde{T_t^{(1)}}$.
The crucial condition for the use of Theorem
\ref{THM:oscintthm} is the boundedness of derivatives of
$\widetilde{a}$ in \eqref{eq:a-tilde}. Here, every
differentiation of $a$ with respect to $y$ introduces
a factor $x_n^{-1}\lambda$ which is then cancelled in view of
assumption \eqref{EQ:assAMP-low-freq1}. It follows that
$\lfloor(n-1)/\gamma\rfloor+1$ $y$-derivatives of $\widetilde{a}$
are bounded, implying the conclusion of
Theorem \ref{THM:oscintthm}. This yields 
estimate \eqref{EQ:est-Tt2} in the way that is
similar to the
proof of the same estimate for $\widetilde{T_t^{(1)}}$.
\end{pf}

We now give an analogue of Proposition \ref{PROP:oscillatory-convex}
without the convexity assumption. Naturally, in this case
we get slower decay rates at infinity.

\begin{prop} \label{PROP:oscillatory-nonconvex}
Let $\gamma_0\in\N$. Let $T_t$ be an operator defined by
\[
T_t f(x)=\int_\Rn e^{i(x\cdot\xi+ t\va(t,\xi))}
\ \br{1-\psi((1+|t|)\xi)}\ b(t,\xi)\ \widehat{f}(\xi) \ \d\xi,
\]
where $\va(t,\xi)$ is real valued, continuous in $t$, 
smooth in
$\xi\in\Rn\backslash 0$, 
homogeneous of order one in $\xi$,
and such that for some $t_0>0$ and $C>0$ we have
\begin{equation}\label{EQ:ass-phi-nc}
C^{-1}|\xi| \leq \varphi(t,\xi) \leq C|\xi| 
\;\;\text{and}
\;\; |\partial_\xi^\alpha\varphi(t,\xi)|\leq C|\xi|^{1-|\alpha|}
\;\;\text{for all}
\;\; t\geq t_0,\; \xi\not=0,
\end{equation} 
and all $\alpha$ such that $|\alpha|\leq \gamma_0+1.$
Assume that the sets
\begin{equation}\label{EQ:ass-sigma-nc}
\Sigma^t=\{\xi\in\Rn\backslash 0:\varphi(t,\xi)=1\}
\end{equation}  
satisfy  
$\lim\sup_{t\to \infty} \gamma_0(\Sigma^t)\leq\gamma_0$ and
$\lim\inf_{t\to\infty} \varkappa_0(\Sigma^t)>0$.  
Let us suppose that the amplitude
$b(t,\xi)$ satisfies 
\begin{equation}\label{EQ:ass-amp-nc}
|\partial_\xi^\alpha b(t,\xi)|\leq C_\alpha |\xi|^{-|\alpha|}
\quad\text{for all}
\quad |\alpha|\leq 1.
\end{equation} 
Let $1< p\leq 2\leq q< \infty$ be such that
$\frac{1}{p}+\frac{1}{q}=1$.
Then for all $t\geq t_0$ we have the estimate
\[
 \Vert T_t f \Vert_{L^q(\Rn)} \leq 
 C t^{-\frac{1}{\gamma_0}\p{\frac{1}{p}-
 \frac{1}{q}}}
 \Vert f \Vert_{\dot{L}^p_{N_p}(\Rn)},
\]
where $N_p=\left(n-\frac{1}{\gamma_0}
\right)
\left(\frac{1}{p}-\frac{1}{q}\right)$. 
\end{prop}

\begin{pf}
Let us show how the proof of Proposition
\ref{PROP:oscillatory-nonconvex} differs from the proof
of Proposition \ref{PROP:oscillatory-convex}. 
First, we need to prove that 
$|I(t,x)|\leq Ct^{-\frac{1}{\gamma_0}}$, 
$t\geq t_0$, for $I_1(t,x)$ as in \eqref{EQ:I-dec}. 
We note that 
$\gamma_0+1\geq 2$ (actually, as we observed
before, we must have $\gamma_0\geq 2$), 
so to prove the estimate for $I_2(t,x)$
we can show that $|I_2(t,x)|\leq Ct^{-1}.$ 
This can be done by integrating by parts with the same
operator $P$ one time, and using \eqref{EQ:ass-amp-nc} instead of
\eqref{EQ:ass-amp}. As for the proof of the estimate
for $I_1(t,x)$, we can reason in the same way as in
Proposition \ref{PROP:oscillatory-convex} to arrive at the
estimate \eqref{EQ:IJ}, i.e.,
\[
|I_1(t,x)|\leq 
C t^{-n+l}2^{j(n-l)}\int_0^\infty 
\abs{J(2^j \lambda,z_t)
\Psi\p{\frac{\lambda}{t }}
\lambda^{n-1-l}} \, d\lambda
\]
with the same operator
$$
J(\lambda,z_t)=\int_U 
e^{i\lambda(-\nabla_y h_t(z_t)\cdot y+ h_t(y)+t x_n^{-1})}
\widetilde{a}(t,x_n,\lambda y,\lambda h_t(y)) 
\kappa_0(t,x,y) \chi(t,y) \, dy.
$$
Now, instead of \eqref{EQ:estJ} we will show that
\begin{equation}\label{EQ:estJ-nc}
\abs{J(\lambda,z_t)}\leq C(1+\lambda)^{-\frac{1}{\gamma_0}},
\quad \lambda>0.
\end{equation}
Then, taking $l=n-\frac{1}{\gamma_0}$, we get the estimate
$|I_1(t,x)|\leq C t^{-\frac{1}{\gamma_0}}$ in the same way
as in estimate \eqref{EQ:I1tx-est2}. Now, estimate
\eqref{EQ:estJ-nc} follows from Theorem \ref{THM:oscintthm}
in the appendix with $N=1$. Indeed, let us write $J(\lambda,z_t)$
in the form \eqref{EQ:gather1}--\eqref{osc-phase} with
phase
$$F(\rho,z_t,\omega)= h_t(\rho\omega+z_t)-h_t(z_t)-
\rho \nabla_z h_t(z_t)\cdot\omega.$$ 
Now, by rotation we may assume that in some direction,
say $e_1=(1,0,\ldots,0)$, we have
by definition of the index $\gamma_0$ that
$$
\gamma_0=\min\{k\in\N: \partial_{\omega_1}^k F(\rho,z_t,\omega)|_
{\omega_1=0}\not=0\}.
$$
Then by taking $N=1$ and $y=\omega_1$ in Theorem 
\ref{THM:oscintthm}, we get the required estimate
\eqref{EQ:estJ-nc}. Conditions 
$\lim\sup_{t\to \infty} \gamma_0(\Sigma^t)\leq\gamma_0$ and
$\lim\inf_{t\to\infty} \varkappa_0(\Sigma^t)>0$
ensure that the dependence on $t$ in this argument is uniform,
so that the constants are also uniform in $t$.
Finally, the estimate for $I_2(t,x)$ differs from that for
$I_1(t,x)$ in exactly the same way as in the proof of
Proposition \ref{PROP:oscillatory-convex}, so we can omit
the repetition of the argument there.
\end{pf}

Suitable $L^p$--$L^q$ estimates for solutions to \eqref{eq:CP} under Assumptions (A1) to (A4) can now be expressed as corollary of these two propositions (depending whether all Fresnel surfaces are convex or not).

\begin{thm}\label{THM:main}
Assume (A1)$_{\ell_1,\ell_2}$, (A2), (A3) and (A4). Assume further that the Fresnel surfaces $\Sigma_k^t$, $k=1,\ldots, m$, are  convex for $t\ge t_0$ for some time $t_0$ and that $\limsup_{t\to\infty} \gamma(\Sigma^t_k)\le \gamma$ with $\liminf_{t\to\infty}\varkappa(\Sigma^t_k)>0$, together with $\lfloor(n-1)/\gamma \rfloor \le \min(\ell_1-1,\ell_2/2-2)$. Then solutions to \eqref{eq:CP} satisfy the $L^p$--$L^q$ estimate
$$
  \| U(t,\cdot)\|_{L^q(\R^n)} \le C_{pq\epsilon} 
  (1+t)^{-\frac{n-1}\gamma\left(\frac1p-\frac1q\right)+\epsilon} 
  \| U_0\|_{W^{p,r_p}(\R^n)} ,
$$
for $p\in(1,2]$, $pq=p+q$, $r_p=n(\frac1p-\frac1q)$, and with 
$\epsilon>0$ arbitrary small if $\nu\in[0,1)$ and for
some $\epsilon>0$ if $\nu=1$.
\end{thm}
\begin{pf}
Due to Theorem~\ref{thm:repSol} solutions to \eqref{eq:CP}  
microlocalised to $Z_{reg}(N,\nu)$ have Fourier integral 
representations of the kind used in 
Proposition~\ref{PROP:oscillatory-convex} after 
multiplication with $t^{-\epsilon}$
(this multiplication makes symbol estimates uniform in $t$ 
as it eliminates the occurring sub-polynomial / polynomial 
terms in the estimates from Theorem~\ref{thm:repSol}). 
It remains to check the assumptions of 
Proposition~\ref{PROP:oscillatory-convex}. 
The phase functions $\phi_j(t,\xi)$ are averages over 
$\lambda_j(t,\xi)$ and therefore homogeneous and bounded 
from below and from above, in view of Propositions
\ref{PROP:characteristics} and \ref{PROP:positivity}. 
Similar for the bounds of their 
derivatives in $\xi$. The amplitudes satisfy the assumptions
 provided $\ell_1$ and $\ell_2$ are big enough, i.e., 
 $\lfloor (n-1)/\gamma\rfloor +1 \le \min(\ell_1,\ell_2/2-1)$,
 combining \eqref{EQ:reg-B-ests} and \eqref{EQ:ass-amp}. 

Finally, combining the microlocal estimate with 
Lemma~\ref{LEM:low-freq} and Sobolev embedding 
theorem to estimate small times concludes the proof.
\end{pf}

We will give one example. If we follow Section~\ref{expl:3} and consider homogeneous 
hyperbolic equations of higher order \eqref{eq:hh-CP},
\begin{equation}
   \D_t^m u + \sum_{j+|\alpha|=m} a_{j,\alpha}(t) \D_t^j \D_x^\alpha u=0,  \quad \D_t^j u (0,\cdot) =u_j, \; j=0,1,\ldots m-1,
\end{equation}  
for $a_{j,\alpha}\in\mathcal T_{\nu}\{0\}$, satisfying the assumption of uniform strict hyperbolicity
in combination with \eqref{EQ:hh-int}, then Theorem~\ref{THM:main} applies and dispersive type estimates depend on geometric properties of the Fresnel surfaces associated to the problem. 
Hence, if  the problem is rotationally invariant all Fresnel surfaces are given by spheres and we obtain $\gamma(\Sigma^t_k) = 2$. The uniformity condition
 $\liminf_{t\to\infty}\varkappa(\Sigma^t_k)>0$ is satisfied if the spheres stay bounded which is 
 the case if $\varphi(t,\xi)$ is uniformly bounded away from zero. In this case the $L^p$--$L^q$ estimate
 \begin{equation}
 \sum_{j=0}^{m-1} \|  |\D|^{m-j-1}  \D_t^j u(t,\cdot) \|_{L^q} \le C_{p,\epsilon} (1+t)^{\epsilon-\frac{n-1}2} \sum_{j=0}^{m-1} \| |\D|^{m-j-1} u_j\|_{W^{p,r_p}}
 \end{equation}
 for  $1<p\le 2$, $pq=p+q$, $r_p = n(\frac1p-\frac1q)$ and
 with arbitrarily small $\epsilon>0$ follows.  If we drop rotational invariance, 
 examples for $\gamma(\Sigma_k^t) \in \{ 2,3, \ldots 2\lfloor m/2 \rfloor\}$ 
 can be constructed in analogy and give corresponding weaker decay rates.

If some of the surfaces fail to be convex, decay rates can be much weaker.

\begin{thm}\label{THM:main-nonconvex}
Assume (A1)$_{1,2}$, (A2), (A3) and (A4), and assume also
that $\limsup_{t\to\infty} \gamma_0(\Sigma^t_k)\le \gamma_0$ 
with $\liminf_{t\to\infty}\varkappa_0(\Sigma^t_k)>0$. 
Then solutions to \eqref{eq:CP} satisfy the $L^p$--$L^q$  estimate
$$
  \| U(t,\cdot)\|_{L^q(\R^n)} \le C_{pq\epsilon} 
  (1+t)^{-\frac{1}{\gamma_0}\left(\frac1p-\frac1q\right)+\epsilon} 
  \| U_0\|_{W^{p,r_p}(\R^n)} 
$$
for $p\in(1,2]$, $pq=p+q$, $r_p=n(\frac1p-\frac1q)$, and with 
$\epsilon>0$ arbitrary small if $\nu\in[0,1)$
and for
some $\epsilon>0$ if $\nu=1$.
\end{thm}

\begin{rem}
In comparison to results on scalar second order equations due to Reissig and co-authors, e.g., \cite{RY98}, \cite{RY00a} or \cite{RS05a}, we observe an $\epsilon$-loss of decay even for the case $\nu=0$.  This is due to the $\xi$-dependence of the term $F_0$ in general.
If $\nu=0$ the choice $\epsilon=0$ can be made in both Theorems \ref{THM:main}
and \ref{THM:main-nonconvex} provided that $\nabla_\xi F_0(t,\xi)=0$.
\end{rem}

\begin{appendix}
\section{The multi-dimensional van der Corput lemma}
We now give the 
multidimensional version of the van der Corput lemma
used in the essential way in
the proof of Proposition \ref{PROP:oscillatory-convex},
as well as in the proof of Proposition \ref{PROP:oscillatory-nonconvex}.

\begin{thm}[\cite{Ru09}]\label{THM:oscintthm}
Consider the oscillatory integral
\begin{equation*}\label{EQ:genoscint}
I(\lambda,\nu)=\int_{\R^{N}}e^{i\lambda
\Phi(x,\nu)}a(x,\nu)\chi(x)\,dx\,,
\end{equation*}
where $N\geq 1$\textup{,} 
and $\nu$ is a parameter. 
Let $\gamma\ge2$ be an integer.
Assume that
\begin{itemize}
\item[(A1)]\label{HYP:mainoscintgbdd} there exists a 
sufficiently small $\delta>0$ such that 
$\chi\in C^\infty_0(B_{\delta/2}(0))$, where
$B_{\delta/2}(0)$ is the ball with radius 
${\delta/2}$ around $0$\textup{;}
\item[(A2)]\label{HYP:mainoscintImPhipos} $\Phi(x,\nu)$ is a
complex valued function such that 
$\textrm{ Im }\Phi(x,\nu)\ge0$ for all
$x\in \supp\chi$ and all parameters $\nu$\textup{;}
\item[(A3)]\label{HYP:mainoscintFconvexfn} for some fixed
$z\in\supp\chi$,  the function
\begin{equation*}
F(\rho,\omega,\nu):=\Phi(z+\rho\omega,\nu), \; |\omega|=1,
\end{equation*}
satisfies the following conditions. Assume that
for each $\mu=(\omega,\nu)$, function $F(\cdot,\mu)$ is of class
$C^{\gamma+1}$ on $\supp\chi$, and
let us write its $\gamma^{\rm th}$ order
Taylor expansion in $\rho$ at $0$ as
\begin{equation*}\label{EQ:Fformwithremainder}
F(\rho,\mu)=\sum_{j=0}^\gamma a_j(\mu)\rho^j + 
R_{\gamma+1}(\rho,\mu)\,,
\end{equation*}
where $R_{\gamma+1}$
is the remainder term. Assume that 
we have
\begin{itemize}
\item[(F1)]\label{ITEM:AssumpF1} $a_0(\mu)=a_1(\mu)=0$ for all
$\mu$\textup{;}
\item[(F2)]\label{ITEM:AssumpF2} there exists a constant $C>0$ such that
$\sum_{j=2}^\gamma\abs{a_j(\mu)}\ge C$ for all
$\mu$\textup{;}
\item[(F3)]\label{ITEM:AssumpF3} for each
$\mu$\textup{,} $\abs{\pa_\rho F(\rho,\mu)}$ is
increasing in $\rho$ for $0<\rho<\delta$\textup{;}
\item[(F4)]\label{ITEM:AssumpF4} for each $k\leq \gamma+1$,
$\pa_\rho^kF(\rho,\mu)$ is bounded uniformly in
$0\leq \rho<\delta$ and $\mu$\textup{;}
\end{itemize}
\item[(A4)]\label{HYP:mainoscintAderivsbdd} for each multi-index $\alpha$
of length $\abs{\alpha}\le \big[\frac{N}{\gamma}\big]+1$\textup{,} there
exists a constant $C_\alpha>0$ such that 
$\abs{\pa_x^\alpha a(x,\nu)}\le
C_\alpha$ for all $x\in \supp\chi$ and all parameters $\nu$.
\end{itemize}
Then there exists a constant $C=C_{N,\gamma}>0$ such that
\[ 
\abs{I(\lambda,\nu)}\le C(1+\lambda)^{-\frac{N}{\gamma}}
\quad\text{for all
}\; \lambda\in[0,\infty)
\textrm{ and all parameters } \nu.
\]
\end{thm}

%
%
%
%

\end{appendix}


\begin{thebibliography}{99}

\bibitem{BL} J.~Bergh, J.~L\"ofstr\"om,
 {\em Interpolation spaces}, Springer, 1976.

\bibitem{Bre75} P. Brenner,
On $L^p$--$L^q$ estimates for the wave equation,
{\em Math. Z.},  {\bf 145}  (1975), 251--254.

\bibitem{dAbbico:2009}
M. d'Abbico, S. Lucente and G. Taglialatela,
$L^p$--$L^q$ estimates for regularly linear hyperbolic systems,
{\em Adv. Diff. Equations}, {\bf 14} (2009), 801--834.\\
Erratum. {\em Adv. Diff. Equations}, {\bf 16} (2011), 199--200.

\bibitem{HW08} F. Hirosawa, J. Wirth,
{$C^m$-theory of damped wave equations with stabilisation}
{\em J. Math. Anal. Appl.}, {\bf 343} (2008), 1022--1035.
%

\bibitem{MR09}
T. Matsuyama, M. Ruzhansky, Time decay for hyperbolic equations 
with homogeneous symbols, 
{\em C. R. Acad. Sci. Paris}, Ser I.  {\bf 347} (2009), 915--919.

\bibitem{MR10} T. Matsuyama, M. Ruzhansky,
Asymptotic integration and dispersion for hyperbolic equations, 
{\em Adv. Diff. Equations}, {\bf 15} (2010), 721--756.

\bibitem{RY98} M. Reissig, K. Yagdjian,
{\em $L_p$--$L_q$ estimates for the solutions of hyperbolic equations of second order
with time-dependent coefficients -- oscillations via growth}, Preprint 98-5, Fakult\"at f\"ur Mathematik und Informatik,
TU Bergakademie Freiberg, 1998.

\bibitem{RY00} M. Reissig, K. Yagdjian, 
The $L^p$--$L^q$ decay estimates for the solutions of 
strictly hyperbolic
equations of second order with increasing in time 
coefficients, {\em Math. Nachr.}, {\bf 214} (2000), 71--104.

\bibitem{RY00a} M. Reissig, K. Yagdjian,
{About the influence of oscillations on Strichartz type decay estimates},
{\em Rend. Sem. Mat. Univ. Pol. Torino}, {\bf 58} (2000), 375--388.

\bibitem{RS05a} M. Reissig, J. Smith,
$L^p$--$L^q$ estimate for wave equations with bounded 
time dependent coefficient,
{\em Hokkaido Math. J.},  {\bf 34}  (2005), 541--586.

\bibitem{Reissig2010} M. Reissig,
Optimality for the asymptotic behaviour of the energy for wave models.
In M. Ruzhansky, J. Wirth (Ed.) {\em Modern Aspects of Partial Differential Equations}, Birkh\"auser, 2011. 


\bibitem{Ru09}
M. Ruzhansky,
Pointwise van der Corput lemma for functions of 
several variables, 
{\em Funct. Anal. Appl.}, 
{\bf 43} (2009), 75--77.

\bibitem{RS05}
M. Ruzhansky, J. Smith,
Global time estimates for solutions to equations of 
dissipative type,
{\em Journ\'ees ``\'Equations aux D\'eriv\'ees Partielles'', 
Exp. No. XII}. Ecole Polytech., Palaiseau, 2005.

\bibitem{RS10}
M. Ruzhansky, J. Smith,
{\em Dispersive and Strichartz estimates for 
hyperbolic equations with constant coefficients}, 
MSJ Memoirs, 22, Mathematical Society of Japan, Tokyo, 2010.

\bibitem{RW08}
M. Ruzhansky, J. Wirth, 
Dispersive estimates for T-dependent hyperbolic systems, 
{\em Rend. Sem. Mat. Univ. Pol. Torino}, {\bf 66} (2008), 339--349.

\bibitem{Str70} R. S. Strichartz,
A priori estimates for the wave equation and some applications,
{\em J. Funct. Anal.}, {\bf 5} (1970), 218--235.

\bibitem{sugi94}
M.~Sugimoto, 
{A priori estimates for higher order hyperbolic equations},
 {\em  Math. Z.}, \textbf{215} (1994), 519--531.

\bibitem{Sug96} M. Sugimoto,
\newblock{Estimates for hyperbolic equations with non-convex characteristics},
\newblock{\em Math. Z.}, {\bf 222} (1996), 521--531.

\bibitem{Wir06} J. Wirth,
Wave equations with time-dependent dissipation I.
 Non-effective dissipation. 
{\em J. Differ. Equations}, {\bf 222} (2006), 487--514.

\bibitem{Wir10} J. Wirth,
Diagonalisation schemes and applications.
{\em Ann. Mat. Pura Appl.}, {\b 189} (2010) 571--590.

\bibitem{Wir10a} J. Wirth,
Energy inequalities and dispersive estimates for wave 
equations with time-dependent coefficients,
{\em Rend. Mat. Univ. Trieste}, {\bf 42 suppl.} (2010) 205--219.

\bibitem{Yagdjian} K. Yagdjian,  {\it The Cauchy problem 
for hyperbolic operators. Multiple characteristics. 
Micro-local approach.} 
Mathematical Topics, 12. Akademie Verlag, Berlin, 1997.  

\end{thebibliography}
\end{document}